\documentclass[10pt,leqno]{amsart} % use larger type; default would be 10pt

\usepackage[english]{babel}
\usepackage{geometry}
\usepackage{mathtools}
\usepackage{amssymb}
\usepackage[hidelinks,final]{hyperref}
\usepackage{graphicx}
\usepackage[dvipsnames]{xcolor}
\usepackage{caption}
\usepackage{float}
\usepackage{subcaption}
\graphicspath{{./Figures/}}
\usepackage[dvipsnames]{pstricks}
\usepackage{pstricks-add}
\usepackage{pst-solides3d}
\usepackage{tikz-cd}
\usepackage{enumerate}
\usepackage{comment}
\usetikzlibrary{babel}
\title{Characterizing maximal varieties via Bredon cohomology}

\author[P.~F.~dos Santos]{Pedro F.~dos Santos}
\email{pedro.f.santos@tecnico.ulisboa.pt}
\address{Departamento de Matemática\\
Instituto Superior Técnico\\
Av. Rovisco Pais\\
1049-001 Lisbon\\
Portugal
}

\author[C.~Florentino]{Carlos Florentino}
\email{caflorentino@fc.ul.pt}
\address{Departamento de Matemática - FCUL\\
Building C6 - 1st floor\\
Universidade de Lisboa\\
Campo Grande, 1749-016 Lisbon\\
Portugal}

\author[J.~Orts]{Javier Orts}
\email{javier.orts@tecnico.ulisboa.pt}
\address{Departamento de Matemática\\
Instituto Superior Técnico\\
Av. Rovisco Pais\\
1049-001 Lisbon\\
Portugal
}

\makeatletter
\numberwithin{equation}{section}
\theoremstyle{plain}
  \newtheorem{thm}{\protect\theoremname}[section]
\theoremstyle{plain}
  \newtheorem{cor}[thm]{\protect\corollaryname}
\theoremstyle{definition}
  \newtheorem{ex}[thm]{\protect\examplename}
  \newtheorem{defi}[thm]{\protect\definitionname}
  
\theoremstyle{plain}
  \newtheorem{prop}[thm]{\protect\propositionname}
\theoremstyle{remark}
  \newtheorem{rem}[thm]{\protect\remarkname}
    \newtheorem{notation}[thm]{\protect\notationname}
\theoremstyle{plain}
  \newtheorem{lem}[thm]{\protect\lemmaname}

\makeatother

  \providecommand{\corollaryname}{Corollary}
  \providecommand{\examplename}{Example}
  \providecommand{\lemmaname}{Lemma}
  \providecommand{\propositionname}{Proposition}
  \providecommand{\remarkname}{Remark}
  \providecommand{\theoremname}{Theorem}
  \providecommand{\definitionname}{Definition}
  \providecommand{\propertyname}{Property}
  \providecommand{\notationname}{Notation}
  \providecommand{\conjecturename}{Conjecture}

\makeatletter
\DeclareFontEncoding{LS1}{}{}
\DeclareFontSubstitution{LS1}{stix}{m}{n}
\DeclareMathAlphabet{\mathscr}{LS1}{stixscr}{m}{n}
\makeatother

\DeclareMathOperator{\id}{id}

\DeclareMathOperator{\Hom}{Hom}

\newcommand{\Z}{\mathbf Z}
\newcommand{\Q}{\mathbf Q}
\newcommand{\R}{\mathbf R}
\newcommand{\C}{\mathbf C}

\newcommand{\F}{\mathbf F}
\newcommand{\RP}{\mathbf{R}P}

\newcommand{\mbP}{\mathbf{P}}

\newcommand{\uFt}{\underline{\F}_2}

\newcommand{\op}{\mathrm{op}}

\newcommand{\M}{\mathbf M}
\newcommand{\Mop}{\M^\star}
\newcommand{\A}{\mathbf A}
\newcommand{\Aop}{\A^{\star}}

\newcommand{\RO}{\mathrm{RO}}

\newcommand{\rk}{\operatorname{Rk}}
\newcommand{\rka}{\operatorname{Rk}_a}
\newcommand{\forget}{\psi}

\newcommand{\leftrarrows}{\mathrel{\raise.75ex\hbox{\oalign{%
  $\scriptstyle\leftarrow$\cr
  \vrule width0pt height.5ex$\hfil\scriptstyle\relbar$\cr}}}}
\newcommand{\lrightarrows}{\mathrel{\raise.75ex\hbox{\oalign{%
  $\scriptstyle\relbar$\hfil\cr
  $\scriptstyle\vrule width0pt height.5ex\smash\rightarrow$\cr}}}}
\newcommand{\Rrelbar}{\mathrel{\raise.75ex\hbox{\oalign{%
  $\scriptstyle\relbar$\cr
  \vrule width0pt height.5ex$\scriptstyle\relbar$}}}}

\makeatletter
\def\leftrightarrowsfill@{\arrowfill@\leftrarrows\Rrelbar\lrightarrows}
\newcommand{\xleftrightarrows}[2][]{\ext@arrow 3399\leftrightarrowsfill@{#1}{#2}}
\makeatother
\newcommand{\hooklongrightarrow}{\lhook\joinrel\longrightarrow}

\DeclareRobustCommand\longtwoheadrightarrow{\relbar\joinrel\twoheadrightarrow}

\usepackage{pst-plot}
\psset{algebraic,arrows=<->}
\definecolor{light-gray}{gray}{0.95}

\renewcommand\qedsymbol{$\blacksquare$}

\begin{document}
\begin{abstract}
We obtain a characterization of Maximal and Galois-Maximal $C_2$-spaces 
(including real algebraic varieties) in terms of $\RO(C_2)$-graded cohomology with coefficients in the constant Mackey functor $\uFt$, using the structure theorem of \cite{clover_may:structure_theorem}.
Other known characterizations, for instance in terms of equivariant Borel cohomology,
are also rederived from this.
%We use an equivariant version of Poincar{\'e} duality from \cite{pedro&paulo:quaternionic_algebraic_cycles} in the the particular case of smooth projective real varieties to deduce further symmetry restrictions for the decomposition of the $\RO(C_2)$-graded cohomology of their complex locus given by the same structure theorem.
For the particular case of a smooth projective real variety $V$, equivariant Poincar\'{e} duality from \cite{pedro&paulo:quaternionic_algebraic_cycles} is used to deduce further symmetry restrictions for the decomposition of the $\RO(C_2)$-graded cohomology of the complex locus $V(\C)$ given by the same structure theorem.
We illustrate this result with some computations, including the $\RO(C_2)$-graded cohomology with $\uFt$ coefficients of real $K3$ surfaces.
\end{abstract}


\maketitle

{\footnotesize\tableofcontents}

%----------------------------------------------------------------------------------------
%	SECTION FILES
%----------------------------------------------------------------------------------------
\section{Introduction}
\label{introduction}

%If $G$ is a $p$-group ($p$ a prime number) and $X$ is a finite $G$-CW complex, Smith theory  provides a relation between the singular cohomology of the space $X$ and that of its fixed points when taking coefficients in the field $\F_p$.

Let $C_2$ denote the group with two elements and let $X$ be a $C_2$-space.
%For a  $X$ be a finite $C_2$-CW-complex\footnote{A finite CW-complex equipped with a $C_2$-action that permutes cells and such that every fixed cell is pointwise fixed.} \textcolor{red}{{\tt Suggestion:} As it is only the intro, we can be sloppy and speak about ``nice'' spaces. We can introduce $C_2$-spaces later on, in Sec. 2}
%
%Let $X$ be a finite CW-complex equipped with an action of $C_2$, the cyclic group of order 2.
%When the action permutes CW-cells \textcolor{red}{(instead, it's cellular)}, we call $X$ a $C_2$-CW complex \textcolor{red}{(fixed cells must be fixed point-wise!)}.
%For such a complex of dimension $n$, 
The well-known Smith--Thom inequality formula \cite{smith:transformations_finite_period}
%,smith:transformations_finite_periodii,smith:fixed-point_theorems} 
%\cite[Theorem 3.3.6]{mangolte:real_algebraic_varieties}
\begin{equation}
\sum_{q=0}^{n} \dim_{\F_2} H_{sing}^q(X^{C_2};\F_2) \quad \leq \quad \sum_{q=0}^{n} \dim_{\F_2} H_{sing}^q(X;\F_2)
\label{eq:ST1}
%\tag{M}
\end{equation}
holds under mild assumptions, thus providing a relation between the singular cohomology of $X$ and that of its fixed points $X^{C_2}$ when taking coefficients in the field $\F_2$.
The spaces for which the equality is attained are called \emph{Maximal spaces}, or \emph{M-spaces} for short.

The natural action of $C_2$ on $H_{sing}^q(X;\F_2)$ turns this cohomology space into a $C_2$-group and allows the following refinement of the Smith--Thom inequality \cite{krasnov:harnack-thom_inequalities}:
\begin{equation}
\sum_{q=0}^{n} \dim_{\F_2} H_{sing}^q(X^{C_2};\F_2) \quad \leq \quad \sum_{q=0}^{n} \dim_{\F_2} H^1(C_2,H_{sing}^q(X;\F_2)). 
\label{eq:ST2}
%\tag{GM}
\end{equation}
Spaces for which this inequality becomes an equality are referred to as \emph{Galois-Maximal spaces} or \emph{GM-spaces} for short.
The notions of M- and GM-spaces originated in real algebraic geometry, where $X$ is the complex locus $V(\C)$ of real algebraic variety $V$, endowed with the analytic topology together with the involution $\sigma$ induced by complex conjugation. % (see Section \ref{subsec:generalities_M_GM}).
%Maximal and Galois-Maximal varieties
They have been the object of study in real algebraic geometry for many years, especially in the context of classification problems (see \cite{silhol:real_algebraic_surfaces}). 

While the concept of Maximal variety dates back to Harnack's work (see \cite{wilson:hilberts_sixteenth_problem} for a survey on Harnack's inequality and Hilbert's XVI problem), Galois-Maximal varieties were introduced in the 1980's by Krasnov in \cite{krasnov:harnack-thom_inequalities} as a generalisation of the former.
Its relation to the classification problem stems from the fact that a real variety $V$ is a GM-variety \emph{iff} the cohomology of $V(\R)$ --- with $\F_2$ coefficients --- is determined by the $C_2$-action on the cohomology of $V(\C)$
%\textcolor{red}{Are we sure about this?} 
(see Section \ref{subsec:generalities_M_GM}). 

Recently there has been renewed interest in the  M-condition in connection with the more restricted notion of \emph{Hodge expressive}
variety (\cite{BruSch21}, \cite{fu23}), focusing on the development of methods to construct new families of M-varieties.

Many of the results on the two maximality conditions for a $C_2$-space $X$ are expressed in terms of the cohomology of the Borel construction $X_{C_2}$.
For example, it is shown in \cite{franz:symmetric_products} that $X$ is an $M$-space \emph{iff} $H^*(X_{C_2};\F_2)$ is a free module over $H^*(BC_2;\F_2).$

In algebraic geometry, $H^*(X_{C_2};\F_2)$ is usually called the $C_2$-equivariant cohomology.
However, in equivariant topology, the $C_2$-equivariant cohomology theory with $\mod 2$ coefficients more commonly used is $\RO(C_2)$-graded equivariant cohomology with coefficients in the constant Mackey functor $\uFt$ (see Section \ref{sec:GM-condition}). 
This theory is less well-known and more complex (as it is bigraded), but it does encode more information and recent developments have made it more accessible to non-experts (\emph{e.g.}, \cite{clover_may:structure_theorem}, \cite{clover_may:freeness_theorem_cohomology}, \cite{Kr09}). 

In this paper, we characterize the maximality conditions in terms of $\RO(C_2)$-graded equivariant cohomology with coefficients in the constant Mackey functor $\uFt$. For simplicity, we call it \emph{bigraded Bredon cohomology}; for a $C_2$-space $X$, we denote it by $H_{C_2}^{*,*}(X;\uFt)$. The bigraded Bredon cohomology ring of a point is denoted $\M_2$.
We refer to the cohomology of $X_{C_2}$ as the \emph{Borel cohomology} of $X$ and denote it by $H^*_{Bor}(X;\F_2)$.

Our main tool is a recent result of Clover May \cite[Thm.~5.1]{clover_may:freeness_theorem_cohomology}. It is a classification for the $\M_2$-modules that can be obtained as the bigraded Bredon cohomology of finite $\C_2$-CW-complexes. 

Using it, we show in Theorem~\ref{thm:main_result} that a finite 
$C_2$-CW-complex $X$ is an M-space \emph{iff} $H_{C_2}^{*,*}(X;\uFt)$ is a free $\M_2$-module, and that $X$ is a GM-space \emph{iff} $H_{C_2}^{*,*}(X;\uFt)$ is the sum of a free $\M_2$-module with a finite number of suspensions of the cohomology of the free orbit $C_2$.

From the bigraded Bredon characterization, we also derive conditions on the Borel cohomology of $X$ for it to be an M- or GM-space.
In the case of M-spaces, this recovers the results of \cite[Prop.~2.1 and Prop.~2.2] {franz:symmetric_products} in the case of a $C_2$-action. 

Finally, in case $X$ satisfies the  bigraded version of Poincar\'{e} duality (see Section \ref{sec:symmetries_poincare}), we show in Theorem~\ref{thm:restrictions_decomposition} that there are symmetries in the decomposition of $H_{C_2}^{*,*}(X;\uFt)$ determined by May's classification result.
We show that certain known conditions implying the GM-condition for smooth real surfaces as in \cite[A1.7]{krasnov:harnack-thom_inequalities} are a consequence of these symmetries (see Prop.~\ref{prop:krasnov}). We also prove a new similar result for real threefolds -- Proposition~\ref{prop:3fold_GM} -- and illustrate it with two examples.

\subsection{Structure of the paper}
\label{subsec:structure}

The structure of this paper is as follows. In
Section~2 we give the definition of Maximal and Galois-Maximal spaces and briefly review the main characterizations used in real algebraic geometry. We also list some examples and constructions that are known to preserve the maximal conditions. In particular, we  describe an important type of M-varieties known as \emph{Hodge expressive varieties}.  

We start Section~3 with an overview of the main properties of bigraded Bredon cohomology needed for our purposes, especially the results of \cite{clover_may:structure_theorem}. 
In Theorem~\ref{thm:main_result} we proceed to characterize Maximal and Galois-Maximal spaces in terms of bigraded Bredon cohomology and in Proposition~\ref{prop:main_result_borel} we give a characterization in terms of Borel cohomology.
In   Corollary~\ref{cor:forgetful_characteristion_M_GM} we show how  other known characterizations follow from ours and give a new characterization GM-spaces in terms of the forgetful map from Borel cohomology.

In Section~4 we review the definition of a Real manifold, Poincar\'{e} duality and the duality relation between bigraded Bredon homology and cohomology established in \cite{clover_may:structure_theorem}, and we prove Theorem~\ref{thm:restrictions_decomposition}. We give several examples involving real algebraic varieties and, in particular, we apply this theorem  to compute the bigraded Bredon cohomology of real $K3$ surfaces. We finish with a few observations concerning Hodge expressive varieties

\subsection{Conventions}
\label{subsec:conventions} In this paper,
 $C_2$ denotes the group with two elements; the irreducible trivial $C_2$-representation is denoted $\mathbf{1}$, $\sigma$ is the $1$-dimensional sign representation and, for $p\geq q$, $\R^{p,q}:=(p-q)\cdot\mathbf{1}\oplus q\cdot\sigma$. The representation sphere of $\R^{p,q}$  is the $C_2$-space $S^{p,q}$ obtained by adding a fixed point at infinity: $S^{p,q}=\R^{p,q}\cup\{\infty\}.$  
The $n$-sphere with the antipodal $C_2$-action is denoted $S_a^n.$ In particular $S^0_a$ is the free $C_2$-orbit.

If $X$ is a $C_2$-space, $X_+$ denotes the based $C_2$-space obtained by adjoining a fixed base point. Given based $C_2$-spaces $(X,x_0)$ and $(Y,y_0)$, their smash product $X\times Y/(x_0\times Y\cup X\times y_0)$ is denoted $X\wedge Y$. For emphasis, the usual (non-equivariant) cohomology of a $C_2$-space $X$ with coefficients in $R$ is denoted $H_{sing}^*(X;R).$

%There appear three different (co)homology theories: singular (co)homology, Borel (co)homology and bigraded Bredon (co)homology.
%The notation is as follows: $H_{sing}^*$, $H_*^{sing}$, $H_{Bor}^*$ and $H_*^{Bor}$ denote singular cohomology and homology and Borel cohomology and homology, respectively; for bigraded (co)homology, this is a bigraded theory, and the notation follows the same of \cite{clover_may:structure_theorem}, \emph{i.e.}~$H_{C_2}^{*,*}$ denotes bigraded cohomology and $H_{*,*}^{C_2}$ denotes bigraded homology.

%The symbols $\M_2$ and $\A_n$ are as in \cite{clover_may:structure_theorem}: $\M_2$ is the $\RO(C_2)$-graded equivariant cohomology ring with coefficients in the constant Mackey functor $\underline{\F_2}$ of a point, and $\A_n$ is $\RO(C_2)$-graded equivariant cohomology group of the antipodal $n$-sphere, $n\geq 0$.

\section*{Acknowledgments}
\label{sec:acknowledgments}
We thank Florent Schaffhauser and Clover May for useful discussions during the preparation of this paper. This work has been partially supported by FCT projects UIDB/04561/2020 and  through CAMGSD, IST-ID, projects UIDB/04459/2020 and UIDP/04459/2020.
The third author wishes to acknowledge the support of FCT through a PhD scholarship.

%\textcolor{red}{I HAVE COMMENTED SOME SECTIONS BECAUSE THE COMPILING TIME EXCEEDS THE LIMIT DUE TO THE GREAT NUMBER OF FIGURES.}

%-----------------------------------------------------------------------------------------------------------------------------------------------
%
%	GENERALITIES ON M- AND GM-SPACES
%
%-----------------------------------------------------------------------------------------------------------------------------------------------
\section{Generalities on maximal spaces and varieties}
\label{subsec:generalities_M_GM}

%Throughout, a CW-complex of finite dimension is assumed to be locally finite 

In this Section, we set up notation and summarize some facts about Maximal and Galois-Maximal spaces. These concepts were first motivated by real algebraic geometry, but they can be defined 
for a large class of spaces which includes finite CW-complexes with an action of $C_2$.
\begin{defi}[G-CW complex]
%\label{ef:G-CW_complex}
Let $G$ be a finite group.
A $G$-CW complex is a CW-complex $X$ on which $G$ acts and satisfies:
\begin{enumerate}

\item The action is cellular, \emph{i.e.}~it sends $n$-cells onto $n$-cells respecting their boundaries;

\item Cells that are fixed under the action are point-wise fixed.

\end{enumerate}
\end{defi}

\begin{ex}
If $V$ is a real algebraic variety (\emph{i.e.}, defined over $\R$) then, by the triangulation results of \cite{Hir75}, its set of complex points $V(\C)$ endowed with the analytic topology and with the $C_2$-action given by the involution $\sigma\colon V(\C)\to V(\C)$ induced from complex conjugation is a $C_2$-CW complex. The set of fixed points $V(\C)^{C_2}$ is the set of real points of $V$, and it is denoted by $V(\R)$. 
\end{ex}
 
\begin{defi}[Real space]
Following \cite{atiyah:k-theory_and_reality}, we call a topological space with an involution a \emph{Real space}. If $V$ is a real algebraic variety, we say that the pair $(V(\C),\sigma)$ is a \emph{Real} algebraic variety (in the sense of Atiyah). 
We call it the Real variety associated to $V$. 
\end{defi}

%-----------------------------------------------------------------------------------------------------------------------------------------------
%
%	M-SPACES AND GM-SPACES
%
%-----------------------------------------------------------------------------------------------------------------------------------------------
\subsection{M-spaces and GM-spaces}

Within the general formalism of $C_2$-spaces, the notions of Maximal and Galois-Maximal spaces come from 
Smith theory (see e.g~ Section 3.2 of \cite{mangolte:real_algebraic_varieties} or 
Chapter III in \cite{bredon:introduction_compact_transformation_groups}), 
and had already been considered by Degtyarev \cite{degtyarev:stiefel_orientations} and Kalinin \cite{kalinin:cohomology_real_algebraic_varieties} before.

\begin{defi}[Maximal spaces]
\label{def:maximal_spaces}
A $C_2$-CW complex $X$ of dimension $n$ is said to be:
\begin{enumerate}[(1)]

\item \emph{Maximal}, abbreviated \emph{M-space}, if the Smith--Thom inequality \eqref{eq:ST1} is an equality:
\begin{equation}
\sum_{q=0}^{n} \dim_{\F_2} H_{sing}^q(X^{C_2};\F_2)  = \sum_{q=0}^{n} \dim_{\F_2} H_{sing}^q(X;\F_2);
\label{eq:M}
\tag{M}
\end{equation}

\item \emph{Galois-Maximal}, abbreviated \emph{GM-space}, if the inequality \eqref{eq:ST2} is an equality:
\begin{equation}
\sum_{q=0}^{n} \dim_{\F_2} H_{sing}^q(X^{C_2};\F_2)  =   \sum_{q=0}^{n} \dim_{\F_2} H^1(C_2,H_{sing}^q(X;\F_2)). 
\label{eq:GM}
\tag{GM}
\end{equation}
\end{enumerate}
\end{defi}

The term $H^1(C_2,A)$ in \eqref{eq:GM} refers to the group cohomology of $C_2$ with coefficients in $A$.
We are interested in the case when $A$ is an $\F_2$-vector space equipped with an involution $\sigma$, for which
\[
H^p(C_2,A)=A^{C_2} \Big/ \{a+\sigma a \mid a\in A\}
\]
(\emph{cf.}~\cite[pp.~55--56]{brown:cohomology_of_groups}).
Hence $\dim_{\F_2} H^1(C_2,A) \leq \dim_{\F_2}A$.
In particular,
\begin{equation*}
%\begin{split}
\sum_{q=0}^{n}\dim_{\F_{2}}H^{1}(C_2,H_{sing}^{q}(X;\F_{2}))
 \leq \sum_{q=0}^{n}\dim_{\F_{2}}H_{sing}^{q}(X;\F_{2}),
%\end{split}
\label{eq:GM-variety_M-variety}
\end{equation*}
(\emph{cf.}~\cite[Lemma 3.6.1]{mangolte:real_algebraic_varieties}) which means that every Maximal space is also Galois-Maximal; the converse does not hold in general.

\begin{defi}[Maximal and Galois-Maximal varieties]
%\label{def:}
A real algebraic variety $V$ is called \emph{Maximal} (resp.~\emph{Galois-mMximal}), if its associated Real variety is Maximal (resp.~Galois-Maximal).
\end{defi}

%The underlying formalism of these notions is the action of the group $C_2$ acting on a topological space (in this particular case the complex locus of the Real variety) by means of the involution $\sigma$.
%Therefore, similar inequalities hold for a much wider class of topological spaces, for instance finite $C_2$-CW complexes.
%It is in this way that one speaks of maximal and Galois-maximal spaces (M- and GM-spaces for short, resp.).

%The precise relation between the two concepts is the following lemma of Krasnov.

%\begin{lem}[Krasnov \cite{krasnov:harnack-thom_inequalities}]
%\label{lem:M-GM}
%Suppose that $X$ is a Galois-Maximal $C_2$-CW-complex.
%Then $X$ is an M-space if and only if $C_2$ acts trivially on the graded 
%$\F_2$-vector space $H_{sing}^*(X;\F_2)$.
%\end{lem}

In real algebraic geometry, GM-spaces are usually characterized in terms of the Leray--Serre
spectral sequence associated to the Borel construction:
if  $X$ is a $G$-space, then 
the \emph{Borel construction} of $X$ is the space
$X_G := (X\times EG)/G$, where $EG\to BG$ is a model for 
the universal $G$-bundle. 
It can be regarded as the total space 
of the fibration
\[
\label{eq:Borel}
X\hooklongrightarrow X_G \longtwoheadrightarrow BG.
\]
We call the singular cohomology of $X_G$ with coefficients in a ring $R$ the \emph{Borel cohomology} of $X$ with coefficients in $R$, with $G$ being understood, and we write:
\[
H_{Bor}^*(X;R) = H_{sing}^*(X_G;R).
\]

The Leray--Serre spectral sequence for the fibration \eqref{eq:Borel} has as its second page
\begin{equation}
\label{eq:spectral-sequence}
E_2^{p,q}=H^p(G,H_{sing}^q(X;\F_2)),
\end{equation}
and converges to the Borel cohomology $H_{Bor}^{p+q}(X;R)=H_{sing}^{p+q}(X_G;R)$.
The inclusion of the fibre $X\hookrightarrow X_G$ provides a map in cohomology $\psi_{Bor}\colon H_{Bor}^*(X; R)\to H_{sing}^*(X;R)$, which we call the \emph{forgetful map}.

We concentrate  on the case $G=C_2$ and $R=\F_2$.

The following Proposition summarizes the characterizations of M-spaces and GM-spaces frequently used in the context of  real algebraic geometry.

%\begin{rem}
%\label{rem:}
%Spaces for which this is onto are called \emph{totally nonhomologous to zero spaces} with respect to the ring $R$ (\emph{cf.}~\cite{borel:seminar_transformation_groups,tom_dieck:transformation_groups}) or \emph{equivariantly formal spaces} (\emph{cf.}~\cite{franz:symmetric_products}).
%For $G=C_2$ and $R=\F_2$ this notion agrees with that of an M-space.
%The reader is again referred to Chapter III in \cite{tom_dieck:transformation_groups}.
%\end{rem}

\begin{prop}[Krasnov \cite{krasnov:harnack-thom_inequalities} and  Franz \cite{franz:symmetric_products}]
\label{prop:M-GM}
Let $X$ be a finite $C_2$-CW complex. Then:
\begin{enumerate}
    \item $X$ is an M-space \emph{iff} the forgetful map
    \[
    H^*_\text{Bor}(X;\F_2)\to H^*_\text{sing}(X;\F_2)
    \]
    is onto;
    \item $X$ is a GM-space \emph{iff}
    the spectral sequence \eqref{eq:spectral-sequence} associated to $X_{C_2}$ 
degenerates at the second page.
\end{enumerate}
In particular, a GM-space $X$ is an M-space \emph{iff}  $C_2$ acts trivially on  $H_{sing}^*(X;\F_2).$
\end{prop}
\begin{proof}
 $ (1)$  is precisely the assertion of \cite[Proposition~2.1]{franz:symmetric_products} in the case $G=C_2$. The two remaining assertions are proved in \cite[Proposition~2.4]{krasnov:harnack-thom_inequalities}.
Although in the original statement $X$ is the  complex locus of a real algebraic variety, the result extends to finite $C_2$-CW complexes, with the same proof.
%
%The proof is essentially the same as that of Proposition~2.4 
%in \cite{krasnov:harnack-thom_inequalities}, since 
%the arguments used can be extended to finite $C_2$-CW-complexes (we refer also to results (1.17), %(4.16) and (4.20) in Chapter III of \cite{tom_dieck:transformation_groups}).
\end{proof}

\subsection{Examples and constructions preserving Maximal and Galois-Maximal varieties}
\label{subsec:maximal_varieties}
We list a few families of real varieties that are known to be Maximal or Galois-Maximal.

\begin{ex}[Maximal and Galois-maximal varieties] 
\phantom{space}\
\label{ex:M_GM_spaces}
\begin{enumerate}
\item Smooth complete curves and Abelian varieties are always GM-varieties, 
provided they have at least one real point \cite{krasnov:harnack-thom_inequalities};
\item Grassmannians (\emph{e.g.}~projective spaces) are Maximal;
\item Smooth complete toric varieties are Maximal \cite{b06}. 
\end{enumerate}
\end{ex}

The following constructions are known to preserve the \eqref{eq:M} condition (resp. \eqref{eq:GM}).

\begin{ex}{Maximal and Galois-Maximal preserving constructions}
\label{ex:constructions-M-GM}
\ \begin{enumerate}
\item Finite products of M-spaces (GM-spaces) are M-spaces (resp. GM-spaces);
\item If $V$ is an M-variety (resp.~GM-variety) and $E\to V$ is a Real vector bundle, then the projective bundle $\mathbf{P}(E)\to V$ is an M-variety (resp.~GM-variety) \cite{fu23};
\item The blow-up of a smooth M-variety along a smooth real M-subvariety is an M-variety \cite{fu23};
\item If $X$ is an M-space then all its symmetric products $\operatorname{SP}_n(X)$ are M-spaces
\cite{franz:symmetric_products}.

\item In a forthcoming paper, we show that all the symmetric products of a GM-space are also GM-spaces.

\end{enumerate}
\end{ex}

\subsection{Hodge expressive varieties}

%By a Real\footnote{Following Atiyah \cite{atiyah:k-theory_and_reality}, we capitalise the letter \textit{R} in the word \textit{Real} when referring to a complex space with an involution.} variety we mean an algebraic variety $V$ over the field of complex numbers that has been equipped with an involution $\sigma:V\to V$ which covers the map $\Spec\C\to\Spec\C$ induced by complex-conjugation (other equivalent definitions of a Real algebraic variety can be found in \cite[Section2.4]{mangolte:real_algebraic_varieties}).
%The fixed points of $V(\C)$ under the action of $\sigma$ are called \emph{real points}, and
%the set of real points, denoted by $V(\R)$, is referred to as the \emph{real locus} of $V$.

%If $V$ is a quasi-projective Real variety, then the pair $(V(\C),V(\R))$ equipped with the Euclidean topology is triangulable (\emph{cf.}~\cite[p.~131]{mangolte:real_algebraic_varieties}), the action of $C_2$ is well-behaved and all the singular cohomology groups are finitely generated. In this situation, Smith theory applies and maximal varieties are defined as before.

In many known examples of smooth M-varieties, the \eqref{eq:M} condition is a consequence
of a more intriguing notion relating Hodge theory to the topology of the set of real points.
We describe it briefly.

Denote by $b_i(X):= \dim_{\F_{2}}H_{sing}^{q}(X;\F_{2})$ 
the $i$-th Betti number of a finite CW-complex $X$ with coefficients in $\F_2$, and by $P_{X}(t)=\sum_{i=0}^d b_i(X)\, t^i$ the corresponding Poincar\'{e} polynomial. Then, for a real variety $V$, the maximal condition can be written using this polynomial as:
%$$
%P_{V(\C)}(t) := \sum_{k=0}^d b_k(V(\C)) \, t^k,
%$$
%where $b_k(V(\C))=\dim H^k(V(\C),\Q)$, since the condition \eqref{eq:M} of 
%Definition \ref{def:maximal_spaces} becomes:
\begin{equation}
\label{eq:M-Poincare}
P_{V(\R)}(1) = P_{V(\C)}(1).
\end{equation}
Now, recall that the rational cohomology of a smooth complex projective variety $V$ of dimension $d$ admits a Hodge decomposition
$$
H^k_{sing}(V(\C);\Q)= \bigoplus_{p+q=k} H^{p,q}(V(\C)),
$$
and this defines the \emph{Hodge polynomial of $V$}:
$$
H_V(u,v) := \sum_{p,q=0}^d h^{p,q}(V(\C)) \, u^p v^q,
$$
where $h^{p,q}(V(\C))=\dim_{\C} H^{p,q}(V(\C))$. Note that, \emph{if $H_{sing}^*(V(\C);\Z)$ is torsion-free,} then
\begin{equation}
\label{eq:torsion-free-hodge-polynomial}
P_{V(\C)}(t) := H_{V(\C)} (t,t).
\end{equation}

\begin{defi}[Hodge-expressive variety]
A real projective variety $V$ is said to be \emph{Hodge-expressive} if $H_{sing}^*(V(\C);\Z)$ is torsion-free and 
\begin{equation}
\label{eq:Hodge-Exp}
P_{V(\R)} (t) = H_{V(\C)} (t,1).
\end{equation}
\end{defi}

\begin{rem}
From \eqref{eq:M-Poincare}, \eqref{eq:torsion-free-hodge-polynomial} and \eqref{eq:Hodge-Exp},
it follows that Hodge expressive varieties are Maximal. Taking the coefficient of $t^i$ in \eqref{eq:Hodge-Exp}, we see that a Hodge expressive variety verifies: 
\begin{equation}
b_i(V(\R)) = \sum_{k=0}^d h^{i,k}(V(\C)).
\label{eq:hodge_expressive_bis}
%\tag{\eqref{eq:Hodge-Exp} bis}
\end{equation}
This is a requirement for a certain compatibility between the topology of $V(\R)$ and the Hodge decomposition that goes beyond the \eqref{eq:M} condition.
\end{rem}

\section{Characterising M- and GM-varieties via equivariant cohomology}
\label{sec:GM-condition}

In this Section, we give a characterization of M-varieties and GM-varieties in terms of two equivariant cohomology theories: $\RO(C_2)$-graded ordinary cohomology and Borel cohomology.
We will start by giving a brief overview of the $\RO(C_2)$-graded theory, and address Borel cohomology in Section \ref{sec:Borel_cohomology}.
To distinguish between the two, we shall refer to the former as \emph{bigraded Bredon cohomology}.

\begin{comment}%Preevious introduction
In this Section, we give a characterization of M-varieties and GM-varieties in terms of bigraded
C2-equivariant cohomology. To distinguish it from other equivariant cohomology theories - like Borel
cohomology - we call it bigraded Bredon cohomology. We start by giving a brief overview of this
cohomology theory.
\end{comment}

\subsection{Preliminaries on bigraded Bredon cohomology}
In \cite{bredon:LMN} Bredon defined an equivariant cohomology theory
$H^n_G(X;M)$ for $G$-spaces, where $G$ is a finite group and $M$
is a contravariant coefficient system. If $M$ is a Mackey functor (see \ref{subsubsection:coeffcients} below), his theory can be extended to an $\RO(G)$-graded theory 
\[ 
\{H^\alpha_G(X;{M}),   \alpha \in \RO(G)\},
\]
called \emph{$\RO(G)$-graded ordinary equivariant cohomology theory},
where $\RO(G)$ denotes the orthogonal representation ring of $G$ --- see \cite{MayJP81}.
%In recent years, the use of $\RO(G)$-graded ordinary equivariant cohomology has become standard in equivariant topology \cite{Kerwaire_invariant}.
 
When $G=C_2$, one has $\RO(C_2) = \mathbb{Z} \cdot\boldsymbol{1} \oplus\mathbb{Z}\cdot \sigma$, where $\boldsymbol{1}$ is the trivial representation and $\sigma$ is the sign representation. In this paper we use $H_{C_2}^{p,q}(X; M)$ to denote $H^{(p-q)\cdot \boldsymbol{1} + q\cdot
\sigma}_{C_2}(X;M)$ and call it \emph{bigraded Bredon cohomology}.
As usual, when $X$ is a based $C_2$-space, we use $\tilde{H}_{C_2}^{p,q}(X;M)$ for the reduced version of the theory. 

Given integers $p\geq q\geq 0$ and a based $C_2$-space $X$, we define 
$\Sigma^{p,q}X:=S^{p,q}\wedge X$, and call it the suspension of $X$ in the direction of the representation sphere $S^{p,q}$.
By analogy with the suspension isomorphism in usual cohomology, bigraded Bredon cohomology satisfies a bigraded suspension axiom:
\[
\label{eq:suspension}
\tilde{H}_{C_2}^{*,*}(\Sigma^{p,q}X;M) =H^{*-p,*-q}_{C_2}(X;M).
\]

\subsubsection{Coefficients} 
\label{subsubsection:coeffcients}
A  Mackey functor $M$ for the group $C_2$ can be described as a diagram of abelian groups of the form,
\[
\begin{tikzcd}
M(C_2)\ar[loop left]{}{t^*}\ar[rr, shift left=1.5ex, "p_*"]&& \ar[ll, "p^*"] M(\bullet)
\end{tikzcd}
\]
where \(C_2\) and \(\bullet=C_2/C_2\) are the types of \(C_2\)-orbits and \(p_*, p^*\) and \(t^*\) are homomorphisms satisfying
\begin{enumerate}
\item $(t^*)^2=\id$;
\item $t^*\circ p^*=p^*$;
\item $p_*\circ t^*=p_*$;
\item $p^*\circ p_*=1+t^*$.
\end{enumerate}

We will use the Mackey functor constant at $\F_2$, which is denoted $\uFt$ and is given by 
\[
\begin{tikzcd}
\F_2\ar[loop left]{}{\id}\ar[rr, shift left=1.5ex, "0"]&& \ar[ll, "\id"] \F_2 .
\end{tikzcd}
\]

\subsubsection{Cohomology of a point: $\M_2$}
The bigraded Bredon cohomology of a point with $\uFt$ coefficients is a bigraded commutative ring  denoted $\M_2$, and described by the following:
As a bigraded abelian group we have $\M_2=\M_2^+\oplus\M_2^-$, with $\M_2^+=\oplus_{p,q\geq 0}\M_2^{p,q}$ and $\M_2^-=\oplus_{p,q<0}\M_2^{p,q}$.
We refer to $\M_2^+$ and $\M_2^-$ as the positive and negative cones of $\M_2$, respectively.
The positive cone is a subring and, as such, it is isomorphic to the polynomial ring $\F_2[\rho,\tau]$, with $\rho\in\M_2^{1,1}$ and $\tau\in \M_2^{0,1}$; the negative cone is more cumbersome to describe. Its 
additive structure is completely determined by the following facts: (1) there is a unique non-zero element $\theta\in \M_2^{0,-2}$; (2) $\rho\theta=0$ and $\tau\theta=0$; (3) all non-zero elements are uniquely divisible by $\rho$ and $\tau$; (4) all elements are obtained from $\theta$ by division with elements of the form $\rho^r\tau^s$.
The multiplicative structure is determined by the fact that it is commutative and that $\theta^2=0$.
Therefore, the product of any two elements of $\M_2^-$ vanishes.
The additive structure of $\M_2$ is depicted in Figure~\ref{fig:M_2}, where every dot indicates a copy of $\F_2$ and the absence of a dot represents the trivial group.    

\subsubsection{Examples of $\M_2$-modules and bigraded rank}
In analogy with other equivariant cohomology theories, bigraded Bredon cohomology with $\underline{\F_2}$ coefficients of a $C_2$-space is always a module over $\M_2$.  
In \cite{clover_may:structure_theorem}, Clover May proved a classification result for the bigraded Bredon cohomology of finite $C_2$-CW complexes, showing that it is a direct sum of only two types of $\M_2$-modules. 
These are described in the next two examples.

If $P$ is a bigraded module over $\M_2$, we denote by $\Sigma^{p,q}P$ the bigraded $\M_2$-module 
with homogenous components given by
$ 
\left(\Sigma^{p,q} P\right)^{u,v} = P^{u-p,v-q}.
$  
Similarly, if $V$ is an $\F_2$-vector space, the $\Z$-graded $\F_2$-space concentrated in degree $p$ is denoted by $\Sigma^pV$.

\begin{figure}[H]
     \centering
     \hfill
     \begin{subfigure}[t]{0.45\textwidth}
         \centering
\begin{pspicture}(-3.5,-3)(3.5,3)
\psset{unit = 0.6}
    \psaxes[linewidth=.4pt,linecolor=gray,showorigin=false,ticksize=2pt,tickcolor=gray,labels=none]{->}(0,0)(-5,-5)(5,5)[$u$,0][$v$,90]
    \psline[linewidth=.7pt,linecolor=black]{->}(0,0)(4.5,4.5)
    \psline[linewidth=.9pt,linecolor=black]{->}(0,0)(0,4.5)
    %\psline[linewidth=.9pt,linecolor=black]{->}(0,-4.5)(0,-2)
    %\psline[linewidth=.7pt,linecolor=black]{->}
    \psline[linewidth=.7pt,linecolor=black,ArrowInside=->]{-}(-2.5,-4.5)(0,-2)
    \psline[linewidth=.7pt,linecolor=black,ArrowInside=->]{-}(0,-4.5)(0,-2)
    
    %\psplot[linecolor=black]{0}{3.3}{x}
%    \rput(*2.5 {x-0.5}){$\times\rho$}
%    \rput(*-0.05 {-999*x}){$\times\tau$}
    \uput[0](1,1){$\rho$}
    \uput[0](-1,1){$\tau$}
    \uput[0](0,-2){$\theta$}
    \uput[0](3.6,3.5){$\times\rho$}
    \uput[0](-1.5,3.5){$\times\tau$}
%    \uput[0](-3.5,-3.5){$\times1/\rho$}
%    \uput[0](-0.1,-3.5){$\times1/\tau$}
    \pscircle[fillstyle=solid,fillcolor=black,dimen=inner](0,0){0.07}
    \pscircle[fillstyle=solid,fillcolor=black,dimen=inner](1,1){0.07}
    \pscircle[fillstyle=solid,fillcolor=black,dimen=inner](2,2){0.07}
    \pscircle[fillstyle=solid,fillcolor=black,dimen=inner](3,3){0.07}
    \pscircle[fillstyle=solid,fillcolor=black,dimen=inner](0,1){0.07}
    \pscircle[fillstyle=solid,fillcolor=black,dimen=inner](0,2){0.07}
    \pscircle[fillstyle=solid,fillcolor=black,dimen=inner](0,3){0.07}
    \pscircle[fillstyle=solid,fillcolor=black,dimen=inner](1,2){0.07}    
    \pscircle[fillstyle=solid,fillcolor=black,dimen=inner](1,3){0.07}
    \pscircle[fillstyle=solid,fillcolor=black,dimen=inner](2,3){0.07}
    \pscircle[fillstyle=solid,fillcolor=black,dimen=inner](0,-2){0.07}
    \pscircle[fillstyle=solid,fillcolor=black,dimen=inner](0,-3){0.07}
    \pscircle[fillstyle=solid,fillcolor=black,dimen=inner](-1,-3){0.07}
    \pscircle[fillstyle=solid,fillcolor=black,dimen=inner](-2,-4){0.07}
    \pscircle[fillstyle=solid,fillcolor=black,dimen=inner](-1,-4){0.07}
    \pscircle[fillstyle=solid,fillcolor=black,dimen=inner](0,-4){0.07}

\end{pspicture}
\caption{\small Representation of the ring $\M_2=\M_2^+\oplus\M_2^-$, as in \cite{clover_may:structure_theorem}.
Multiplication by $\rho$ occurs diagonally, while multiplication by $\tau$ is along the vertical axis.}
         \label{fig:M_2}
     \end{subfigure}
     \hfill
     \begin{subfigure}[t]{0.45\textwidth}
         \centering

\begin{pspicture}(-3,-3)(3,3)
\psset{unit = 0.5}
    \psaxes[linewidth=.4pt,linecolor=gray,showorigin=false,ticks=none,labels=none]{->}(-3.5,-3.5)(-6,-5)(6,5)[$u$,0][$v$,90]
    \psline[linewidth=.7pt,linecolor=black]{->}(0,0)(4.5,4.5)
    \psline[linewidth=.85pt,linecolor=black]{->}(0,0)(0,4.5)
    \psline[linewidth=.7pt,linecolor=black,ArrowInside=->]{-}(-2.5,-4.5)(0,-2)
    \psline[linewidth=.7pt,linecolor=black,ArrowInside=->]{-}(0,-4.5)(0,-2)
    
    %\psplot[linecolor=black]{0}{3.3}{x}
%    \rput(*2.5 {x-0.5}){$\times\rho$}
%    \rput(*-0.05 {-999*x}){$\times\tau$}
    \uput[0](1,1){$\rho$}
    \uput[0](-1,1){$\tau$}
    \uput[0](0,-2){$\theta$}
    \uput[0](3.6,3.5){$\times\rho$}
    \uput[0](-1.5,3.5){$\times\tau$}
%    \uput[0](-3.5,-3){$\times1/\rho$}
%    \uput[0](0,-4){$\times1/\tau$}
    \uput[0](0,0){$(p,q)$}

    \pscircle[fillstyle=solid,fillcolor=black,dimen=inner](0,0){0.07}
    \pscircle[fillstyle=solid,fillcolor=black,dimen=inner](1,1){0.07}
    \pscircle[fillstyle=solid,fillcolor=black,dimen=inner](2,2){0.07}
    \pscircle[fillstyle=solid,fillcolor=black,dimen=inner](3,3){0.07}
    \pscircle[fillstyle=solid,fillcolor=black,dimen=inner](0,1){0.07}
    \pscircle[fillstyle=solid,fillcolor=black,dimen=inner](0,2){0.07}
    \pscircle[fillstyle=solid,fillcolor=black,dimen=inner](0,3){0.07}
    \pscircle[fillstyle=solid,fillcolor=black,dimen=inner](1,2){0.07}    
    \pscircle[fillstyle=solid,fillcolor=black,dimen=inner](1,3){0.07}
    \pscircle[fillstyle=solid,fillcolor=black,dimen=inner](2,3){0.07}
    
    \pscircle[fillstyle=solid,fillcolor=black,dimen=inner](0,-2){0.07}
    \pscircle[fillstyle=solid,fillcolor=black,dimen=inner](-1,-3){0.07}
    \pscircle[fillstyle=solid,fillcolor=black,dimen=inner](-2,-4){0.07}
    \pscircle[fillstyle=solid,fillcolor=black,dimen=inner](0,-3){0.07}
    \pscircle[fillstyle=solid,fillcolor=black,dimen=inner](0,-4){0.07}
    \pscircle[fillstyle=solid,fillcolor=black,dimen=inner](-1,-4){0.07}

\end{pspicture}
         \caption{Picture of $\Sigma^{p,q}\M_2$ for $(p,q)$ a real representation.
         This module can be realized as the reduced cohomology of the representation sphere $S^{p,q}$.}

         \label{fig:M_2shifted}
     \end{subfigure}
     \hfill
     \caption{Representation of the modules $\M_2$ and $\Sigma^{p,q}\M_2$.}
        \label{fig:M_2&A_n}
\end{figure}
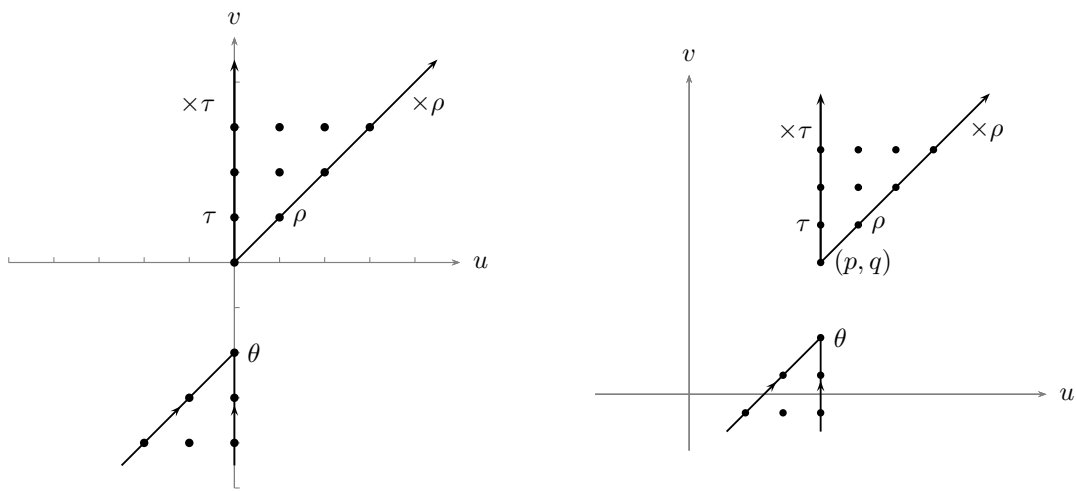

\begin{ex}[The modules $\Sigma^{p,q}\M_2$]
If $P$ is the reduced cohomology of a based $C_2$-space $X$ then, by the bigraded suspension axiom \eqref{eq:suspension},
$\Sigma^{p,q}P$ is the reduced cohomology of the suspension $\Sigma^{p,q}X:$
\[
\widetilde{H}^{*,*}(\Sigma^{p,q}X;\uFt) =\widetilde{H}^{*-p,*-q}(X;\uFt) = (\Sigma^{p,q}P)^{*,*},
\]
which justifies the notation $\Sigma^{p,q}P$ for $\M_2$-modules.
In particular, $\Sigma^{p,q}\M_2=\widetilde{H}^{*,*}(S^{p,q};\uFt)$  is a free $\M_2$-module with a single generator
in bidegree $(p,q)$, as depicted in Figure~\ref{fig:M_2shifted}.
\end{ex}

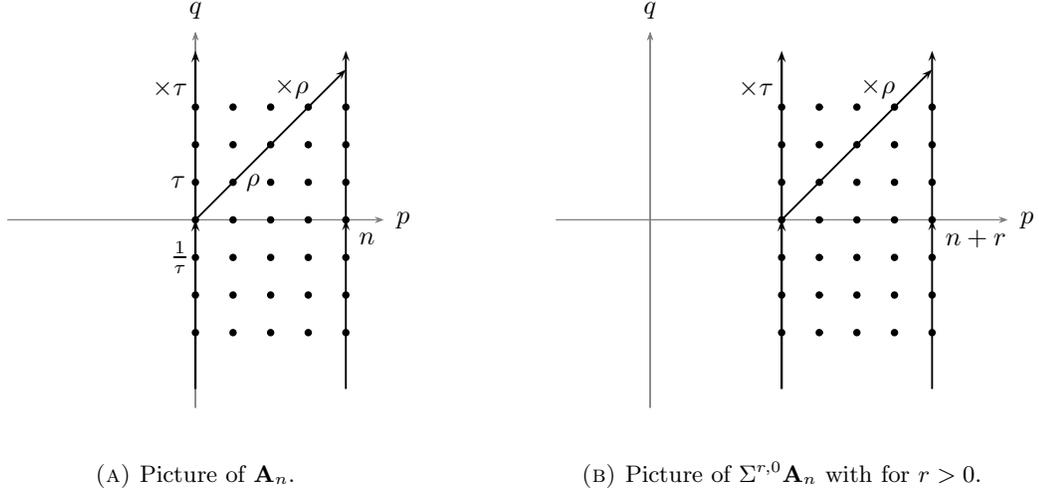
\begin{figure}[H]
     \centering

     \hfill
     \begin{subfigure}[t]{0.45\textwidth}
    
    \centering

\begin{pspicture}(-3,-3)(3,3)
\psset{unit = 0.5}
    \psaxes[linewidth=.4pt,linecolor=gray,showorigin=false,ticks=none,labels=none]{->}(0,0)(-5,-5)(5,5)[$p$,0][$q$,90]
    
    \psline[linewidth=.7pt,linecolor=black]{->}(0,0)(4,4)
    \psline[linewidth=.85pt,linecolor=black]{->}(0,0)(0,4.5)
    \psline[linewidth=.85pt,linecolor=black]{->}(0,-4.5)(0,0)
    \psline[linewidth=.7pt,linecolor=black]{->}(4,-4.5)(4,0)
    \psline[linewidth=.7pt,linecolor=black]{->}(4,0)(4,4.5)

    %\psplot[linecolor=black]{0}{3.3}{x}
    %\rput(*2.5 {x-0.5}){$\times\rho$}
    %\rput(*-0.05 {-999*x}){$\times\tau$}
    \uput[0](1,1){$\rho$}
    \uput[0](-1,1){$\tau$}
    \uput[0](-1,-1){$\frac{1}{\tau}$}
    \uput[0](4,-0.5){$n$}
    
    \uput[0](1.75,3.5){$\times\rho$}
    \uput[0](-1.5,3.5){$\times\tau$}
%    \uput[0](-2.25,-3.5){$\times1/\tau$}
    
    \pscircle[fillstyle=solid,fillcolor=black,dimen=inner](0,0){0.07}
    \pscircle[fillstyle=solid,fillcolor=black,dimen=inner](1,0){0.07}
    \pscircle[fillstyle=solid,fillcolor=black,dimen=inner](2,0){0.07}
    \pscircle[fillstyle=solid,fillcolor=black,dimen=inner](3,0){0.07}
    \pscircle[fillstyle=solid,fillcolor=black,dimen=inner](4,0){0.07}
    
    \pscircle[fillstyle=solid,fillcolor=black,dimen=inner](0,1){0.07}
    \pscircle[fillstyle=solid,fillcolor=black,dimen=inner](1,1){0.07}
    \pscircle[fillstyle=solid,fillcolor=black,dimen=inner](2,1){0.07}
    \pscircle[fillstyle=solid,fillcolor=black,dimen=inner](3,1){0.07}
    \pscircle[fillstyle=solid,fillcolor=black,dimen=inner](4,1){0.07}
    \pscircle[fillstyle=solid,fillcolor=black,dimen=inner](0,2){0.07}
    \pscircle[fillstyle=solid,fillcolor=black,dimen=inner](1,2){0.07}
    \pscircle[fillstyle=solid,fillcolor=black,dimen=inner](2,2){0.07}
    \pscircle[fillstyle=solid,fillcolor=black,dimen=inner](3,2){0.07}
    \pscircle[fillstyle=solid,fillcolor=black,dimen=inner](4,2){0.07}
    \pscircle[fillstyle=solid,fillcolor=black,dimen=inner](0,3){0.07}
    \pscircle[fillstyle=solid,fillcolor=black,dimen=inner](1,3){0.07}
    \pscircle[fillstyle=solid,fillcolor=black,dimen=inner](2,3){0.07}
    \pscircle[fillstyle=solid,fillcolor=black,dimen=inner](3,3){0.07}
    \pscircle[fillstyle=solid,fillcolor=black,dimen=inner](4,3){0.07}
    \pscircle[fillstyle=solid,fillcolor=black,dimen=inner](0,-1){0.07}
    \pscircle[fillstyle=solid,fillcolor=black,dimen=inner](1,-1){0.07}
    \pscircle[fillstyle=solid,fillcolor=black,dimen=inner](2,-1){0.07}
    \pscircle[fillstyle=solid,fillcolor=black,dimen=inner](3,-1){0.07}
    \pscircle[fillstyle=solid,fillcolor=black,dimen=inner](4,-1){0.07}
    \pscircle[fillstyle=solid,fillcolor=black,dimen=inner](0,-2){0.07}
    \pscircle[fillstyle=solid,fillcolor=black,dimen=inner](1,-2){0.07}
    \pscircle[fillstyle=solid,fillcolor=black,dimen=inner](2,-2){0.07}
    \pscircle[fillstyle=solid,fillcolor=black,dimen=inner](3,-2){0.07}
    \pscircle[fillstyle=solid,fillcolor=black,dimen=inner](4,-2){0.07}
    \pscircle[fillstyle=solid,fillcolor=black,dimen=inner](0,-3){0.07}
    \pscircle[fillstyle=solid,fillcolor=black,dimen=inner](1,-3){0.07}
    \pscircle[fillstyle=solid,fillcolor=black,dimen=inner](2,-3){0.07}
    \pscircle[fillstyle=solid,fillcolor=black,dimen=inner](3,-3){0.07}
    \pscircle[fillstyle=solid,fillcolor=black,dimen=inner](4,-3){0.07}

\end{pspicture}
\caption{\small Picture of $\A_n$.}
         \label{fig:A_4}
     \end{subfigure}
     \hfill
     \begin{subfigure}[t]{0.45\textwidth}
         \centering

\begin{pspicture}(-3,-3)(3,3)
\psset{unit = 0.5}
    \psaxes[linewidth=.4pt,linecolor=gray,showorigin=false,ticks=none,labels=none]{->}(-3.5,0)(-6,-5)(6,5)[$p$,0][$q$,90]

    \psline[linewidth=.7pt,linecolor=black]{->}(0,0)(4,4)
    \psline[linewidth=.85pt,linecolor=black]{->}(0,0)(0,4.5)
    \psline[linewidth=.85pt,linecolor=black]{->}(0,-4.5)(0,0)
    \psline[linewidth=.7pt,linecolor=black]{->}(4,-4.5)(4,0)
    \psline[linewidth=.7pt,linecolor=black]{->}(4,0)(4,4.5)
    
    %\psplot[linecolor=black]{0}{3.3}{x}
    %\rput(*2.5 {x-0.5}){$\times\rho$}
    %\rput(*-0.05 {-999*x}){$\times\tau$}
     %\uput[0](1,1){$\rho$}
    %\uput[0](-1,1){$\tau$}
    %\uput[0](-2,-1){$\tau^{-1}$}
    \uput[0](4,-0.5){$n+r$}
    
    \uput[0](1.75,3.5){$\times\rho$}
    \uput[0](-1.5,3.5){$\times\tau$}
%    \uput[0](-2.25,-3.5){$\frac{1}{\tau}$}
    
    \pscircle[fillstyle=solid,fillcolor=black,dimen=inner](0,0){0.07}
    \pscircle[fillstyle=solid,fillcolor=black,dimen=inner](1,0){0.07}
    \pscircle[fillstyle=solid,fillcolor=black,dimen=inner](2,0){0.07}
    \pscircle[fillstyle=solid,fillcolor=black,dimen=inner](3,0){0.07}
    \pscircle[fillstyle=solid,fillcolor=black,dimen=inner](4,0){0.07}
    
    \pscircle[fillstyle=solid,fillcolor=black,dimen=inner](0,1){0.07}
    \pscircle[fillstyle=solid,fillcolor=black,dimen=inner](1,1){0.07}
    \pscircle[fillstyle=solid,fillcolor=black,dimen=inner](2,1){0.07}
    \pscircle[fillstyle=solid,fillcolor=black,dimen=inner](3,1){0.07}
    \pscircle[fillstyle=solid,fillcolor=black,dimen=inner](4,1){0.07}
    \pscircle[fillstyle=solid,fillcolor=black,dimen=inner](0,2){0.07}
    \pscircle[fillstyle=solid,fillcolor=black,dimen=inner](1,2){0.07}
    \pscircle[fillstyle=solid,fillcolor=black,dimen=inner](2,2){0.07}
    \pscircle[fillstyle=solid,fillcolor=black,dimen=inner](3,2){0.07}
    \pscircle[fillstyle=solid,fillcolor=black,dimen=inner](4,2){0.07}
    \pscircle[fillstyle=solid,fillcolor=black,dimen=inner](0,3){0.07}
    \pscircle[fillstyle=solid,fillcolor=black,dimen=inner](1,3){0.07}
    \pscircle[fillstyle=solid,fillcolor=black,dimen=inner](2,3){0.07}
    \pscircle[fillstyle=solid,fillcolor=black,dimen=inner](3,3){0.07}
    \pscircle[fillstyle=solid,fillcolor=black,dimen=inner](4,3){0.07}
    \pscircle[fillstyle=solid,fillcolor=black,dimen=inner](0,-1){0.07}
    \pscircle[fillstyle=solid,fillcolor=black,dimen=inner](1,-1){0.07}
    \pscircle[fillstyle=solid,fillcolor=black,dimen=inner](2,-1){0.07}
    \pscircle[fillstyle=solid,fillcolor=black,dimen=inner](3,-1){0.07}
    \pscircle[fillstyle=solid,fillcolor=black,dimen=inner](4,-1){0.07}
    \pscircle[fillstyle=solid,fillcolor=black,dimen=inner](0,-2){0.07}
    \pscircle[fillstyle=solid,fillcolor=black,dimen=inner](1,-2){0.07}
    \pscircle[fillstyle=solid,fillcolor=black,dimen=inner](2,-2){0.07}
    \pscircle[fillstyle=solid,fillcolor=black,dimen=inner](3,-2){0.07}
    \pscircle[fillstyle=solid,fillcolor=black,dimen=inner](4,-2){0.07}
    \pscircle[fillstyle=solid,fillcolor=black,dimen=inner](0,-3){0.07}
    \pscircle[fillstyle=solid,fillcolor=black,dimen=inner](1,-3){0.07}
    \pscircle[fillstyle=solid,fillcolor=black,dimen=inner](2,-3){0.07}
    \pscircle[fillstyle=solid,fillcolor=black,dimen=inner](3,-3){0.07}
    \pscircle[fillstyle=solid,fillcolor=black,dimen=inner](4,-3){0.07}
\end{pspicture}
         \caption{\small Picture of $\Sigma^{r,0}\A_n$ with for $r> 0$.}
         \label{fig:A_4shifted}
     \end{subfigure}
     \hfill

\caption{\small As abelian group $\A_n\cong \F_2[\tau,\tau^{-1},\rho]/(\rho^{n+1})$.
%When regarded as a module over $\M_2$, multiplication by $\theta$ annihilates any element of $\A_n$.
}

\label{fig:A_n}
\end{figure}

\begin{ex}[The cohomology of antipodal spheres]
Denote by $\A_n$ the bigraded Bredon cohomology ring of $S^n_a$  (see Section~\ref{subsec:conventions}): $\A^{*,*}_n=H_{\C_2}^{*,*}({S^n_a}_+;\uFt)$. 
In particular, $\A_0$ is the bigraded cohomology of the free $C_2$ orbit.
As an $\M_2$-algebra, $\A_n$  can be described as follows \cite{clover_may:structure_theorem}: $\A_n\cong \F_2[\tau,\tau^{-1},\rho]/(\rho^{n+1})$.
Note that, as multiplication by the degree $(0,1)$ element $\tau$ is invertible, $\A_n$ is $(0,1)$-periodic and $\theta$ annihilates any element of $\A_n$ because $\tau\theta=0$.
It also follows that, for all $p,q$, the $\M_2$-modules $\Sigma^{p,q}\A_n$ and $\Sigma^{p,0}\A_n$ are isomorphic.
The latter appears represented in Figure \ref{fig:A_4shifted}.
\end{ex}

It will be convenient to introduce two \emph{rank functions} on all those $\M_2$-modules which are direct sums of
suspensions of $\M_2$ and $\A_n$.

\begin{defi} 
\label{def:bigraded_ranks}
Let $I$ and $J$ be finite indexing sets, and let $N$ be an $\M_2$-module satisfying a decomposition as follows:
\begin{equation}
\label{eq:standard-sum}
N \cong \bigoplus_{i\in I} \Sigma^{p_i,q_i}\M_2 \;\oplus \;\bigoplus_{j\in J} \Sigma^{r_j,0}\A_{n_j}
\end{equation}
Define
$\rk(N), \rka(N)$ $\colon\Z^2\to\Z_{\geq 0}$ by
\begin{equation}
\rk^{p,q}(N) := |\{i\in I: (p_i,q_i)=(p,q) \}| \quad\text{and}\quad \rka^{s,t}(N):= |\{ j\in J: (r_j, n_j)=(s,t) \} |.
\label{eq:symmetry_restrictions}
\end{equation}
It is clear that  $\rk(N)$ and $\rka(N)$ determine the isomorphism type of $N.$
\end{defi}

\begin{rem}
\label{rem:bigraded_ranks} 
A priori it is not obvious that $\rk(N)$ and $\rka(N)$ are functions of the module $N$ and not of the particular decomposition used in their definitions.
In the case of $\rk(N)$, we can make this clear by expressing it alternatively as
\[
\rk^{p,q}(N) = \dim_{\F_2}(N/IN)^{p,q}
\]
where $I\subset\M_2$ denotes the kernel of the projection $\M_2\to\F_2$. This function was first introduced in \cite{Dug15} for finitely generated free $\M_2$-modules. It is called the \emph{bigraded rank} of $N$. 

Regarding the definition of $\rka(N)$, note that, if $N$ decomposes as in \eqref{eq:standard-sum}, then, by setting $\M_2^{\rk{N}}:=$ $\oplus_{(p,q)}\M_2^{\rk(N)(p,q)}$, we have
\begin{align*}
N/\M_2^{\rk{N}}
&\cong 
\bigoplus_{j\in J} \Sigma^{r_j,0}\A_{n_j}\\
&\cong \bigoplus_{j\in J}\F_2[\tau,\tau^{-1},\rho]/(\rho^{n_j+1}).
\end{align*}
It now follows from the classification theorem for finitely generated modules over the PID $\F_2[\tau,\tau^{-1},\rho]$ that $\rka^{r,s}(N)$ is determined by the isomorphism type of $N$: $\rka^{r,s}(N)$ is the number of $r$-shifted copies of $\F_2[\tau,\tau^{-1},\rho]/(\rho)^{s+1}$ in the primary decomposition of $N$.
We call $\rka(N)$ the \emph{bigraded $a$-rank} of $N$.
\end{rem}

%\begin{rem}
%%\label{rem:}
%If $X$ has a fixed point then there exists at least one index $i\in I$ such that $(p_i,q_i)=(0,0)$, for
%%
%\[
%H_{C_2}^{*,*}(X;\underline{\F_2}) = \M_2 \oplus \tilde{H}_{C_2}^{*,*}(X;\underline{\F_2}).
%\]
%%
%%, say $x_0$, then the inclusion $x_0\hookrightarrow X$ and the projection $X\rightarrow x_0$ provide maps in equivariant cohomology whose composition
%%%
%%\[
%%H_{C_2}^{*,*}(x_0;\underline{\F_2}) \longrightarrow H_{C_2}^{*,*}(X;\underline{\F_2}) \longrightarrow %H_{C_2}^{*,*}(x_0;\underline{\F_2}) 
%%\]
%%%
%%must be the identity.
%%This shows that when $X^{C_2}$ is non-empty,.
%\end{rem}

\subsubsection{Forgetful map to singular cohomology}
An important property of bigraded Bredon cohomology that will be used throughout is its relation to (non-equivariant) singular  cohomology: for a \(C_2\)-space \(X\) there is a natural isomorphism  
\begin{equation}
\label{eq:singular.equivariant.iso}
H^{p,q}_{C_2}(X\times C_2;\uFt)\cong H_{sing}^p(X;\F_2)
\end{equation}
and, under this isomorphism, the equivariant map \(C_2\to C_2/C_2\) induces a natural map
\begin{equation}
\label{eq:forgetful}
\forget\colon H^{p,q}_{C_2}(X;\uFt) \to  H^{p,q}_{C_2}(X\times C_2;\uFt)\cong H_{sing}^p(X;\F_2),
\end{equation}
which we call the \emph{forgetful map}. It is a map of cohomology theories.

\begin{ex}[Forgetful map on representation and antipodal spheres]
On $\M_2$, the forgetful map is completely determined by $\psi(\tau)=1$ and $\psi(\rho)=0.$ The fact that $\psi$ commutes with suspensions 
yields  $\psi\bigl(\Sigma^{p,q}\M_2\bigr)=\Sigma^p\F_2.$
In the same way, from $\A_n=\F_2[\tau,\tau^{-1},\rho]/(\rho^{n+1})$, it follows that $\psi(\Sigma^{r,0}\A_n)=\Sigma^r\F_2.$ In particular, note that $\psi$ is not surjective; it is concentrated in degree $r.$
\end{ex}

\subsubsection{$\rho$-localisation and fixed points}
For a $C_2$-space $X$, the singular cohomology of $X^{C_2}$ can be
computed from the bigraded Bredon cohomology by a localisation formula as follows:
\begin{lem}[Clover May, \cite{clover_may:structure_theorem}]
For a finite \(C_2\)-CW complex, we have
\begin{equation}
\rho^{-1}H_{C_2}^{*,*}(X;\underline{\F_2}) \cong H_{sing}^*(X^{C_2};\F_2)\otimes_{\F_2} \rho^{-1}\M_2,
\label{eq:rho_localisation}
\end{equation}
where elements of degree $k$ in singular cohomology are considered to have bidegree $(k,0).$ 
\end{lem}

\begin{ex}[Fixed points of representation spheres]
\label{ex:fixed_pts_rep_spheres}
Since \(\Sigma^{p,q}\M_2\)
is  a free \(\M_2\)-module with one generator in degree \((p,q)\), we have 
\[
\rho^{-1}\Sigma^{p,q}\M_2\cong \Sigma^{p,q}\rho^{-1}\M_2 \cong \Sigma^{p-q,0}\rho^{-1}\M_2
= \tilde{H}_{sing}^*(S^{p-q,0};\F_2)\otimes_{\F_2}\rho^{-1}\M_2.
\]
\end{ex}

%-----------------------------------------------------------------------------------------------------------------------------------------------
%
%	THE STRUCTURE THEOREM AND MAIN CHARACTERISATION
%
%-----------------------------------------------------------------------------------------------------------------------------------------------

\subsection{The structure theorem and main characterisation}
\label{subsec:structure_theorem}

In \cite{clover_may:structure_theorem}, Clover May provided a complete characterization of the bigraded Bredon cohomology of a finite $C_2$-CW complex with coefficients in the constant Mackey functor $\underline{\F_2}$ which only involves $\M_2$-modules of the form $\Sigma^{p,q}\M_2$ and $\Sigma^{r,0}\A_n$, for $p\geq q, r, n\in\Z_{\geq 0}$.

\begin{thm}[Structure theorem, \cite{clover_may:structure_theorem}]
\label{thm:structure_theorem}
For a finite $C_2$-CW complex X, there is a decomposition of the bigraded Bredon equivariant cohomology of $X$ with constant coefficients $\underline{\F_2}$ as
\begin{equation}
H_{C_2}^{*,*}(X;\underline{\F_2}) \, =\, \bigoplus_{i\in I} \Sigma^{p_i,q_i}\M_2 \; \oplus \; \bigoplus_{j\in J} \Sigma^{r_j,0}\A_{n_j}
\label{eq:structure_theorem}
\end{equation}
as a module over $\M_2$.
Furthermore, $p_i\geq q_i\geq 0$ and $r_j,n_j\geq 0$ for all $i\in I$ and $j\in J$.
\end{thm}

\begin{notation}
Often we will need to isolate the terms in \eqref{eq:structure_theorem} which are suspensions of \(\A_0\). For that purpose, we set
\[
J_0 := \{j\in J: n_j=0\} \quad \text{and} \quad J_+ := \{j\in J: n_j>0\}.
\]
\end{notation}

\begin{ex}
For the projective line $\mbP^1$, the Real space $\mbP^1(\C)$ is the representation sphere $S^{2,1}$; hence
\[
H_{C_2}^{*,*}(\mbP^1(\C);\uFt) =\M_2\oplus\Sigma^{2,1}\M_2.
\]
More generally, in the case of  the projective $n$-space, we have
\begin{equation}
\label{eq:decomp_Pn}
H_{C_2}^{*,*}(\mbP^n(\C);\uFt)= \bigoplus_{0\leq i \leq n }\Sigma^{2i,i}\M_2.
\end{equation}    
\end{ex}

\begin{ex}
\label{ex:torus}
Let $\Lambda\subset \C$ be the lattice  $\Z+\sqrt{-1}\cdot\Z$ and let $E=\C/\Lambda$ denote the corresponding elliptic  curve.
The action of $C_2$ on \(E(\C)\) is induced by complex conjugation on \(\C\). As a \(C_2\)-space \(E(\C)\) is isomorphic to \(S^{1,0}\times S^{1,1}\) and so
%In the geometric realisation, this action is a reflection with respect to the horizontal plane.
%The two circles that result as the intersection of the torus with the horizontal plane remain fixed under this action, providing a factor $\Sigma^{1,0}\M_2$ in bigraded cohomology.
%They are the real locus of $T$.
%Poincaré duality then implies the existence of a factor $\Sigma^{1,1}\M_2$, which can be associated to either of the representation spheres $S^{1,1}$.
the bigraded Bredon cohomology is
\[
H_{C_2}^{*,*}(E(\C);\underline{\mathbf{F}_2}) = \M_2 \oplus \Sigma^{1,0}\M_2 \oplus \Sigma^{1,1}\M_2 \oplus \Sigma^{2,1}\M_2.
\]
\end{ex}

\begin{ex}
A Severi--Brauer variety over a field $K$ is a smooth projective variety such that its extension of scalars to the algebraic closure $\overline K$ is isomorphic to a projective space (over $\overline K$). For $K=\R,$ the non-trivial Severi--Brauer curve is the conic $\operatorname{SB}_1\subset \mathbf{P}^2$ described by 
\[
x^2+y^2+z^2=0.
\]
One can check that the associated Real variety $(\operatorname{SB}_1(\C),\sigma)$ is isomorphic as a $C_2$-space to the antipodal sphere $S^2_a.$ It follows that 
\label{ex:brauer-severi}
\[
H_{C_2}^{*,*}(\operatorname{SB}_1(\C);\uFt) = \A_2.
\]
\end{ex}

Theorem~\ref{thm:structure_theorem} implies that the isomorphism type of the bigraded cohomology of $X$ as an $\M_2$-module is completely determined by the bigraded rank and the $a$-rank of $H_{C_2}^{*,*}(X;\underline{\F_2})$, as in 
Definition \ref{def:bigraded_ranks}. For convenience, we introduce a two-variable polynomial codifying the information contained in the bigraded rank function.

\begin{defi}
\label{def:bigraded_rank_polynomial}
For a finite $C_2$-CW complex, define 
$ 
r^{p,q}(X) := \rk^{p,q}H_{C_2}^{*,*}(X;\underline{\F_2})
$ 
and
$ 
R_X (u,v) := \sum_{p\geq q\geq 0} r^{p,q}(X) \, u^p \, v^q .
$
\end{defi}

%\textcolor{red}{{\tt Notation:} I have changed the notation of the bigraded rank functions to lowercase, in resemblance to the Hodge numbers.}

In order to apply the structure theorem to express the conditions \eqref{eq:M} and \eqref{eq:GM} in terms of the graded Bredon cohomology of the space we need to use May's localisation formula \eqref{eq:rho_localisation} to read from \eqref{eq:structure_theorem} the required information about the singular cohomology of \(X\) and \(X^{C_2}\).

It turns out that the bigraded rank polynomial of $X$ determines the Poincaré polynomial 
of the fixed point locus $X^{C_2}$.
\begin{thm}
\label{thm:cohomology_fixed_points}
Let $X$ be a finite $C_2$-CW complex with equivariant cohomology given by equation \eqref{eq:structure_theorem}. Then
%\begin{enumerate}
%\item There exists an isomorphism %of $C_2$-groups
%
\begin{equation*}
H_{sing}^*(X^{C_2};\F_2) \cong \bigoplus_{i\in I} \tilde{H}_{sing}^*(S^{p_i-q_i};\F_2).
%\oplus \left(\bigoplus_{j\in J} \tilde{H}_{sing}^*(S^{r_j}\wedge (S^{n_j})_+;\F_2)\right)
\label{eq:singular_cohomology}
\end{equation*}
Hence, for each $r\geq 0$, the $k$-th Betti number of $X^{C_2}$ is given by:
\[
b_k(X^{C_2}) \, = \, \dim H_{sing}^k(X^{C_2};\F_2) \, = \, 
|\{i\in I\mid p_i-q_i = k\}| \, = \,  \sum_{m \geq 0} r^{k+m,m}(X)
\]
or, equivalently,
\[
P_{X^{C_2}} (t) \, = \, R_X(t, \, t^{-1} ).
\]
\end{thm}

\begin{proof}
For the second type of  summand in  decomposition 
\eqref{eq:structure_theorem}, \(\Sigma^{r_j,0}\A_{n_j}\),  all elements are \(\rho\)-torsion hence 
their \(\rho\)-localization is trivial. 
From this computation, and from Example~\ref{ex:fixed_pts_rep_spheres}, the first assertion follows. The second assertion is a direct consequence of the first, by using the identity
\[
\sum_{\substack{i\in I: \\ p_i-q_i=k}} r^{p_i,q_i}(X) \, = \, \sum_{m \geq 0} r^{k+m,m}(X).
\]
Finally, an easy computation gives:
\[
R_X(t, t^{-1}) = \sum_{p\geq q\geq 0} r^{p,q}(X) \, t^{p-q} 
= \sum_{k,q \geq 0} r^{k+q,q}(X) \, t^k = 
\sum_{k \geq 0} b_k(X) \, t^k = P_{X^{C_2}} (t),
\]
as desired.
\end{proof}

Theorem \ref{thm:cohomology_fixed_points} implies the following characterization of Hodge expressive varieties.

\begin{cor}
Let $V$ be a real variety which is smooth and projective and whose integer cohomology has no torsion. 
Then, $V$ is Hodge expressive if and only if 
\begin{equation}
\label{eq:Hodge-rank}
H_{V(\C)}(t, 1) = R_{V(\C)}( t, \, t^{-1} ).
\end{equation}
\end{cor}
\begin{proof}
From Theorem \ref{thm:cohomology_fixed_points}, applied to $X=V(\C)$ and $X^{C_2}=V(\R)$, we have:
\[
R_{V(\C)}( t,  \, t^{-1}  ) = P_{V(\R)}( t).
\]
The result now follows from Equation \eqref{eq:Hodge-Exp}. 
\end{proof}

\begin{rem}
\label{rem:Hodge-birank}
Taking coefficients of the polynomials, identity \eqref{eq:Hodge-rank} is equivalent to:
\[
\sum_{q \geq 0} h^{p,q}(V(\C)) \, =\,  \sum_{q \geq 0} r^{p+q,q}(V(\C)).
\]
Hence, the equality $r^{p+q,q}(V(\C))=h^{p,q}(V(\C))$ for all $p,q\geq 0$ is sufficient for Hodge expressivity. All the examples of Hodge expressive varieties in Section \ref{subsec:applications_poincare_symmetries} verify this equality.
So, it is an interesting problem to find an example of Hodge expressive variety which does not satisfy
the above relation between the Hodge numbers and the bigraded rank.
\end{rem}

\subsubsection{Mackey-valued cohomology functor}
As we have seen, bigraded Bredon cohomology and singular cohomology are related by the forgetful map.
Actually, there is a way to package both theories into a single functor by enhancing  \(X\mapsto H^{*,*}_{C_2}(X;\uFt)\)  to a Mackey-functor-valued theory \(\ X\mapsto \underline{H_{C_2}^{*,*}}(X)\) that can be described by the following diagram:
\[
\begin{tikzcd}
%\ar[loop left,"t^*"] 
\ar[loop left]{}{t^*}
\underline{H_{C_2}^{*,*}}(X)(C_2)\hspace{-1cm}&:= H_{C_2}^{*,*}(X\times C_2;\uFt)
\ar[rr, shift left=1.5ex, "p_*"]&& \ar[ll, "p^*"] 
H_{C_2}^{*,*}(X;\uFt) &\hspace{-1cm}=:
\underline{H_{C_2}^{*,*}}(X)(\bullet).
\end{tikzcd}
\]
Note  that, in particular, by \eqref{eq:singular_cohomology} the value of this functor at $C_2$ is identified with singular cohomology.
Under this identification  \(p^*\) is the forgetful map of \eqref{eq:forgetful} and \(t^*\) is the map induced by the \(C_2\)-action on $X$.  The existence of \(p_*\) satisfying the relations in the Mackey functor definition is an important  feature of equivariant cohomology; it is known as the \emph{transfer} corresponding to the cover \(C_2\to\bullet\); see \cite{MayJP96}.
%\footnote{It can be constructed using the equivariant suspension axiom and an  appropriate equivariant map \(S^V\to S^V\wedge {C_2}_+,\) for some \(C_2\)-representation \(V;\) see \cite{MayJP96}.} 

\begin{ex}Denoting the trivial $C_2$ orbit by $\bullet,$ we have $\underline{H_{C_2}^{0,0}}(\bullet)=\uFt.$
\end{ex}

We will not compute this more structured functor in the examples we present, but we do use the fact that May's  structure theorem for bigraded cohomology applies to it.

\bigskip

In order to obtain the singular cohomology of $X$ from decomposition \eqref{eq:structure_theorem}, we need to use the enhanced functor \(\ X\mapsto \underline{H_{C_2}^{*,*}}(X)\). 
As explained in \cite{clover_may:structure_theorem}, the decomposition of \eqref{eq:structure_theorem} is valid for this this Mackey valued functor.

\begin{thm}[Mackey-valued Structure theorem, \cite{clover_may:structure_theorem}] 
Under the same conditions of Theorem \ref{thm:structure_theorem}, we have an isomorphism of Mackey functors
\label{thm:Mackey.struture.thm}
\begin{equation}
\underline{H_{C_2}^{*,*}}(X) \,\cong\,
\bigoplus_{i\in I} \underline{\tilde{H}_{C_2}^{*,*}}(S^{p_i,q_i}) \; \oplus \;
\bigoplus_{j\in J} \underline{\tilde{H}_{C_2}^{*,*}}\left(S^{r_j,0}\wedge {S^{n_j}_a}_{\! +}\right).
\label{eq:Mackey_structure_theorem}
\end{equation}
\end{thm}

\begin{cor}
\label{cor:singular_cohomology}
Let $X$ be a finite $C_2$-CW complex, with graded Bredon cohomology as in Theorem \ref{thm:structure_theorem}. 
Then there is an isomorphism of graded abelian groups with a \(C_2\)-action
\begin{equation}
H_{sing}^*(X;\F_2) \,\cong\, \bigoplus_{i\in I} \Sigma^{p_i}\F_2 \;\oplus \;
\bigoplus_{j\in J_0} \Sigma^{r_j}\F_2[C_2] \;\oplus\;
\bigoplus_{j\in J_+}(\Sigma^{r_j}\F_2 \;\oplus\;
\Sigma^{r_j+n_j}\F_2),
\label{eq:singular_cohomology}
\end{equation}
where
\[
J_0 := \{j\in J: n_j=0\} \quad \text{and} \quad J_+ := \{j\in J: n_j>0\}.
\]
\end{cor}

\begin{proof}
This follows from Theorem \ref{thm:Mackey.struture.thm}, since  the value of the functor at \(C_2\) can be identified with singular cohomology, by \eqref{eq:singular.equivariant.iso}.
\end{proof}

We can now we  state our characterization result.
\begin{thm}
Let $X$ be a finite $C_2$-CW complex under the same conditions of Theorem \ref{thm:structure_theorem}.
\label{thm:main_result}
%Let $X$ be a finite $C_2$-CW complex with bigraded Bredon cohomology decomposition as follows
%\begin{equation}
%H_{C_2}^{*,*}(X;\underline{\F_2}) = \left(\bigoplus_{i\in I} \Sigma^{p_i,q_i}\M_2\right) \oplus %\left(\bigoplus_{j\in J} \Sigma^{r_j,0}\A_0\right) \oplus \left(\bigoplus_{k\in K} \Sigma^{r_k,0}\A_{n_k}\right) 
%\label{eq:strcuture_theorem_bis}
%\end{equation}
Then:
%: with bigraded equivariant cohomology
%%
%\begin{equation*}
%\begin{split}
%H_{C_2}^{*,*}(X;\underline{F_2}) = \left(\bigoplus_{i\in I} \tilde{H}_{C_2}^{*,*}\left(S^{p_i,q_i};\underline{\F_2}\right) \right) & \oplus \left(\bigoplus_{j\in J}\tilde{H}_{C_2}^{*,*}\left(S^{r_j}\wedge (S_a^0)_+;\underline{\F_2}\right)\right)\\
%& \oplus\left(\bigoplus_{k\in K}\tilde{H}_{C_2}^{*,*}\left(S^{r_k}\wedge(S_a^{n_k})_+;\underline{\F_2}\right)\right),
%%\label{eq:strcuture_theorem_bis}
%\end{split}
%\end{equation*}
%%
%where the contributions of $n$-dimensional spheres in the structure theorem are split into those with $n=0$ and $n>0$ for convinience and all $I,J,K$, $p_i,q_i,\dots$ satisfy the same conditions that in Theorem \ref{thm:structure_theorem}.
%Then:

%\renewcommand{\labelenumi}{(\arabic{enumi})}
\begin{enumerate}

\item $X$ is an M-space \emph{iff} it has free bigraded cohomology, 
\emph{i.e.}, \emph{iff} $J=\varnothing$.
\item $X$ is a GM-space \emph{iff} the only non-free summands in decomposition 
\eqref{eq:structure_theorem} are of the form $\Sigma^{r,0}\A_0$, \emph{i.e.}, \emph{iff}   $J_+=\varnothing$.
\end{enumerate}
\end{thm}

\begin{proof}
%By Theorem~\ref{thm:Mackey.struture.thm},
%the singular cohomology of $X$ is
%
%\[
%H_{sing}^*(X;\F_2) \cong_{C_2} \left(\bigoplus_{i\in I} \tilde{H}_{sing}^*(S^{p_i};\F_2)\right) \oplus %\left(\bigoplus_{j\in J} H_{sing}^{*+r_j}(S^{0};\F_2)\right)\oplus \left(\bigoplus_{k\in K}H_{sing}^{*+r_k}%(S^{n_k};\F_2)\right).
%\]
%

%Let $|I|$ denote the cardinality of the indexing set $I$.
%Notice that it is equivalent to the number of (representation) spheres in the cohomology of $X$.
%Similarly, denote by $|J|$ and $|K|$ the number of antipodal $0$-spheres and of antipodal $n$-spheres ($n>0$), respectively.
%
By Corollary~\ref{cor:singular_cohomology},
\[
\sum_{q=0}^{\dim X} \dim_{\F_2} H_{sing}^q(X;\F_2) = |I| + 2|J_0| + 2|J_+|
\]
and, by Theorem~\ref{thm:cohomology_fixed_points}
\[
\sum_{q=0}^{\dim X^{C_2}} \dim_{\F_2} H_{sing}^q(X^{C_2};\F_2) = |I|.
\]
Hence $X$ satisfies \eqref{eq:M} if and only if $J=\varnothing$, which is equivalent to the equivariant cohomology being free.

As for  the GM-property, by \eqref{eq:singular_cohomology}, we have
\[
\sum_{q=0}^{\dim X} \dim_{\F_2} H^1(C_2,H_{sing}^q(X;\F_2)) = |I| + 2|J_+|,
\]
since  $\F_2[C_2]$ has no higher $C_2$-cohomology and so the suspensions of antipodal $0$-spheres do not contribute to this sum.
The result follows.
\end{proof}

We can also characterize the M- and GM-conditions in terms of the image of the forgetful map to singular cohomology. 

\begin{cor}
\label{cor:forgetful-charaterisations_M_GM_conditions}
Let $X$ be a finite $C_2$-CW complex, and let $W\subset H_{sing}^*(X;\F_2)$ be the image of the 
forgetful map $\psi:H_{C_2}^{*,*}(X;\uFt)\to H_{sing}^*(X;\F_2)$.
\begin{enumerate}
    \item $X$ is an M-space \emph{iff} $\psi$ is onto, \emph{i.e.}, $W = H_{sing}^*(X;\F_2)$;
    \item $X$ is an GM-space \emph{iff} $W$ is the fixed point space $H_{sing}^*(X;\F_2)^{C_2}$.
\end{enumerate}
\end{cor}
\begin{proof}
This is clear from Theorem~\ref{thm:main_result}, display \eqref{eq:singular_cohomology} and the evaluation of  forgetful maps on $\Sigma^{p,q}\M_2$ and $\Sigma^{n,0}\A_n.$
\end{proof}

\begin{rem}
\label{rem:rederive_Kransnov}
Note that the above corollary provides an alternative proof of Krasnov's result on the relation between 
M- and GM-spaces stated in Proposition \ref{prop:M-GM}.
Indeed, if $X$ is a GM-space and the $C_2$-action on $H_{sing}^*(X;\F_2)$ is trivial, then the fixed point space is the whole vector space, so 
$W=H_{sing}^*(X;\F_2)$, which means that $X$ is an M-space by the corollary.
\end{rem}

%-----------------------------------------------------------------------------------------------------------------------------------------------
%
%	THE STRUCTURE THEOREM AND MAIN CHARACTERISATION
%
%-----------------------------------------------------------------------------------------------------------------------------------------------

\subsection{Characterization via Borel cohomology}
\label{sec:Borel_cohomology}
The equivariant cohomology theory more commonly used in real algebraic geometry is Borel cohomology.
In this section, we characterize the \eqref{eq:M} and \eqref{eq:GM} conditions in terms of this theory, 
and provide several examples.

Borel cohomology is related to bigraded Bredon cohomology as follows \cite{MayJP96}:
\[
H^{*}_{Bor}(X;\F_2) = H^{*,0}_{C_2}(X\times_{C_2} EC_2;\uFt),
\]
where $EC_2$ denoted the total space for the principal universal $C_2$-bundle.

In general, it  encodes less information than its bigraded Bredon counterpart (see \cite[Ch.IV, \S2]{MayJP96}). 
For finite $C_2$-CW complexes, using Theorem~\ref{thm:structure_theorem}, May obtained a structure theorem for Borel cohomology by $\tau$-localisation, as follows.

\begin{lem}[$\tau$-localisation, \cite{clover_may:structure_theorem}]
\label{lem:tau_localisation}
For any finite $C_2$-CW complex $X$, identifying $z$ with $\rho/\tau$, the following is a natural isomorphism of $\F_2[z]$-modules:
\[
\left[\tau^{-1}H_{C_2}^{*,*}(X;\underline{\F_2})\right]^{*,0} \cong H^*_{Bor}(X;\F_2).
\]
\end{lem}

\begin{ex}[Borel cohomology of a point and of representation spheres]
Noting that $\tau^{-1}\M_2=\F_2[\tau,\tau^{-1},\rho]$, we have 
$  
 H_{Bor}^*(\bullet;\F_2) = \F_2[z],
$    
 where $z=\rho\tau^{-1},$ and
\[ 
 \widetilde{H}_{Bor}^*(S^{p,q};\F_2) = \Sigma^p\F_2[z].
 \]  
Similarly, since $\tau$ is invertible in $\A_n,$ we have 
\[
\widetilde{H}_{Bor}^*(S^{r,0}\wedge {S^n_a}_+;\F_2)=\A_n^{*,0}=\Sigma^r\F_2[z]/(z^{n+1}).
\]
\end{ex}

The existence of a \emph{forgetful map}
\[
\psi_{Bor}\colon H^*_{Bor}(X;\F_2) \to H^*_{sing}(X;\F_2)
\]
is a well known and commonly used feature of Borel cohomology. It can be shown that, under the isomorphism  of Lemma~\ref{lem:tau_localisation}, $\psi_{Bor}$ is the map on the $\tau$-localization induced from $\psi$ (note that $\psi(\tau)=1$). 
In particular, $\psi_{Bor}\bigl(\Sigma^r\F_2[z]/(z^{n+1})\bigr)=\Sigma^r\F_2$
and $\psi_{Bor}\bigl(\Sigma^p\F_2[z]\bigr)=\Sigma^p\F_2$.

\begin{cor}[May {\cite[\S6.16]{clover_may:structure_theorem}}]
\label{cor:structure_corollary}
Let $X$ is a finite $C_2$-CW complex whose bigraded Bredon cohomology decomposes as in \eqref{eq:structure_theorem}:
\begin{equation*}
H_{C_2}^{*,*}(X;\underline{\F_2}) \, =\,  \bigoplus_{i\in I} \Sigma^{p_i,q_i}\M_2 \; \oplus \; \bigoplus_{j\in J} \Sigma^{r_j,0}\A_{n_j},
\end{equation*}
where $p_i\geq q_i\geq 0$ and $r_j,n_j\geq 0$ for all $i\in I$ and $j\in J$.
Then $H_{Bor}^*(X;\F_2)$  has the following decomposition over  $H_{Bor}^*(\bullet;\F_2)=\F_2[z]$:
\begin{equation}
\label{eq:Borel_Stucture}
H^*_{Bor}(X;\F_2) \, =\, \bigoplus_{i\in I} \Sigma^{p_i}\F_2[z] \; \oplus \; \bigoplus_{j\in J}\Sigma^{r_j}\F_2[z]/(z^{n_j+1}).
\end{equation}
\end{cor}

\begin{rem}
This result  shows exactly what information is lost when going from bigraded Bredon cohomology to Borel cohomology.
As explained in Subsection \ref{subsec:structure_theorem}, the bigraded Bredon cohomology of $X$ is determined by the lists of pairs $(p_i,q_i)$, $(r_j,n_j)$, counted with multiplicity. In the language of Definition~\ref{def:bigraded_ranks}, it is determined by  $\rk  H^{*,*}_{C_2}(X;\uFt)$ and $\rka  H^{*,*}_{C_2}(X;\uFt)$.
By Corollary~\ref{cor:structure_corollary}, its Borel cohomology is determined by 
$\rka  H^{*,*}_{C_2}(X;\uFt)$ and
the lists of $p_i$, counted with multiplicity. The information about the $q_i$'s is precisely what is lost in the process. 

%To characterize the information lost in this process we can use the following 
%alternative description of the bigraded rank function using  an increasing filtration on $H_{sing}^*(X;\F_2)$: for each $p,q\in\Z_{\geq 0}$, setting
%\[
%F^q H^{p,*}_{C_2}(X;\uFt) := \bigcup_{q'\leq q}\psi \bigl(H^{p,q'}_{C_2}(X;\uFt)\bigr).
%\]
%we have
%
%\[
%\rk^{p,q} H^{*,*}_{C_2}(X;\uFt) = \dim_{\F_2}\frac{F^q H^{p,*}_{C_2}(X;\uFt)}{F^{q-1} H^{p,*}_{C_2}(X;\uFt)}.
%\]
%
%Dropping the information about the $q_i's$ amounts to keeping only the information about the total dimensional jumps in the filtration on $H_{sing}^p(X;\F_2)$, for each $p.$ 
\end{rem}

Applying  Corollary~\ref{cor:structure_corollary} and Theorem~\ref{thm:main_result} we can now characterize the spaces that
satisfy the \eqref{eq:M} and \eqref{eq:GM} conditions in terms of their Borel cohomology.

\begin{prop}
\label{prop:main_result_borel}
Let $X$ be a finite $C_2$-complex.% and assume that its $C_2$-graded cohomology is exactly as in \eqref{eq:structure_theorem}.
Then:
\begin{enumerate}%[ref=(\arabic*)]
\item $X$ satisfies the \eqref{eq:M} condition \emph{iff} its Borel cohomology is free and finitely generated as a (graded) module over $\F_2[z]$.
\label{item:maximal_spaces_borel}

\item  $X$ satisfies the \eqref{eq:GM} condition \emph{iff} its Borel cohomology has a decomposition over $\F_2[z]$ of the following type:
\begin{equation}
H^*_{Bor}(X;\F_2) \, =\, \bigoplus_{i\in I}\Sigma^{p_i}\F_2[z]\; \oplus \; \bigoplus_{j\in J} \Sigma^{r_j}\F_2 .
\label{eq:borel_cohomologyGM}
\end{equation}
\end{enumerate}
\end{prop}
\begin{proof}
This follows immediately from  Theorem~\ref{thm:main_result} and display \eqref{eq:Borel_Stucture}.
\end{proof}

Statement (\ref{item:maximal_spaces_borel}) in the above proposition is a particular case of \cite[(4.16) of III]{tom_dieck:transformation_groups}) for $G=C_2$.
It follows from \eqref{eq:Borel_Stucture} that 
\[
\psi_{Bor}\bigl( H^p_{Bor}(X;\F_2)\bigr) \, =\, \bigoplus_{i\in i}\Sigma^{p_i}\F_2 \;\oplus\; \bigoplus_{j\in J_0}\Sigma^{r_j}\F_2[C_2]^{C_2}\!.
\]
Hence, it exists a characterization of the M- and GM-conditions in terms of the image of $\psi_{Bor}$.

\begin{cor}
\label{cor:forgetful_characteristion_M_GM}
Let $X$ be a finite $C_2$-CW complex. Then:
\begin{enumerate}
    \item $X$ is an M-space \emph{iff} the \emph{forgetful map}  $\psi_{Bor}$  is onto;
    \item $X$ is an GM-space \emph{iff} the image  of the \emph{forgetful map}  $\psi_{Bor}$ is $H_{sing}^*(X;\F_2)^{C_2}.$
\end{enumerate}
\end{cor}

%% Symmetries and restrictions imposed by Poincaré duality

\section{Symmetries and restrictions imposed by Poincar\'{e} duality}
\label{sec:symmetries_poincare}

\subsection{Homology and Poincar\'{e} duality}
Bigraded Bredon cohomology has a homological counterpart. It is a bigraded theory 
$X\mapsto H_{p,q}^{C_2}(X;\uFt)$ satisfying the covariant versions of the properties satisfied by bigraded Bredon cohomology. There is a class of $C_2$-spaces satisfying a  bigraded version of Poincar\'{e} duality that relates these two  bigraded invariants. We refer to them as \emph{Real manifolds (in the sense of Atiyah)} \cite[Section 1]{pedro&paulo:quaternionic_algebraic_cycles}.

\begin{defi}[Real manifold]
A \emph{Real manifold} of dimension $n$, also called an \emph{$n$-Real manifold}, is a $C_2$-manifold $M$ such that its tangent bundle
$TM$ has a Real $n$-bundle structure. That is, $TM$ has a complex vector bundle structure of dimension $n$ and a complex anti-linear involution covering the involution on $M$. In particular, $M$ is a real manifold of dimension $2n$.
\end{defi}

\begin{ex}[Real varieties as Real manifolds]
For a smooth real $n$-variety $V$, the set of complex points $V(\C)$ equipped with the analytic topology and the involution induced by complex conjugation is a $n$-Real manifold.
\end{ex}

If  $M$ is a compact  Real $n$-manifold, then the Poincaré duality theorem for bigraded Bredon cohomology \cite[\S1.6]{pedro&paulo:quaternionic_algebraic_cycles} states that there are natural isomorphisms of abelian groups
\[
H_{p,q}^{C_2}(M;\uFt) \cong H^{2n-p,n-q}_{C_2}(M;\uFt)
\]
(\emph{cf.}~Proposition 4.14 in \textit{loc.~cit.}).%\cite[Prop.~1.14]{pedro&paulo:quaternionic_algebraic_cycles}).
%Actually, as in the non-equivariant case,  Bredon homology is a module of cohomology 

It is therefore natural to apply a homological version of the structure Theorem~\ref{thm:structure_theorem} in conjunction with Poincar\'{e} duality to obtain further restrictions on the bigraded Bredon cohomology of smooth Real manifolds.

%As an application of these restrictions, we show in the end of this section how to rederive a result of Krasnov about real GM-surfaces --- Proposition~\ref{prop:krasnov} below.

We will start by explaining how the homological version of  Theorem~\ref{thm:structure_theorem} is established in 
\cite{clover_may:structure_theorem}. 
There, it is shown that for a finite $C_2$-CW complex $X$, its bigraded Bredon homology $H_{*,*}^{C_2}(X;\uFt)$ is the dual of $H^{*,*}_{C_2}(X;\uFt)$ in the appropriate sense: for an $\M_2$-module $N,$ its dual $N^\star$ is the $\M_2$-module defined by
\begin{equation}
\label{eq:def_dual}
N^\star_{p,q} := \Hom_{\M_2}\left(N,\Sigma^{p,q}\M_2\right),
\end{equation}
and then it is shown that
\begin{equation}
H_{p,q}^{C_2}(X;\underline{\F_2}) \cong H^{*,*}_{C_2}(X;\underline{\F_2})^\star_{p,q}.
\label{eq:homology_dual_cohomology}
\end{equation}

\begin{rem}
From  \eqref{eq:def_dual}, it is clear that $N^\star$ is naturally an $\M_2$-module, but not  a \emph{bigraded $\M_2$-module}, as $\M_2$ acts with reverse grading: if  $\mu\in\M_2^{p,q}$ and $v\in N^\star_{p',q'}$ then $\mu\cdot v\in N^\star_{p'-p,q'-q}$.
This is to be  expected since, as in the non-equivariant case, there is a cap product making cohomology act on homology with reverse grading.  
\end{rem}

It is easy to check that the operation $\star$ commutes with suspensions and direct sums, hence  it follows that if $X$ is a finite $C_2$-CW complex with cohomology given by \eqref{eq:structure_theorem}, then
\begin{equation}
\label{eq:homology_structure}
H^{C_2}_{*,*}(X;\uFt) \,=\, \bigoplus_{i\in I} \Sigma^{p_i,q_i}\M_2^\star \;\oplus\; \bigoplus_{j\in J} \Sigma^{r_j,0}\A_{n_j}^{\!\star}.
\end{equation}

To describe $\Mop_2$ and $\Aop_n$ the following definition is useful.

\begin{defi}
Given a commutative bigraded ring $R$ we define its \emph{opposite ring}, $R^\op$, as the graded ring obtained by setting $R^\op_{p,q}:=R^{-p,-q}$, with the same operations as $R$.
Similarly, if $N$ is a graded $R$-module, we define $N^\op$ by setting $N^\op_{p,q}:=N^{-p,-q}$, with the same operations. In particular, $N^\op$ is a bigraded $R^\op$-module. 
\end{defi}

\begin{ex}[$\Mop_2$ and $\Aop_n$]
\label{ex:M2_An_duals}
From the definition, it is easy to check that $\Mop_2=\M_2^{\op}$ and $\Aop_n=\Sigma^{n,0}\A_n^{\op}$, as $\M_2$-modules.
Both are depicted in the following figure (Fig.~\ref{fig:duals}):

\begin{figure}[H]
     \centering
     \hfill
     \begin{subfigure}[t]{0.45\textwidth}
         \centering
\begin{pspicture}(-3.5,-3)(3.5,3)
\psset{unit = 0.5}
    \psaxes[linewidth=.4pt,linecolor=gray,showorigin=false,ticks=none,labels=none]{->}(0,0)(-5,-5)(5,5)[$p$,0][$q$,90]
%    \psline[linewidth=.7pt,linecolor=black]{->}(0,2)(2.5,4.5)
%    \psline[linewidth=.9pt,linecolor=black]{->}(0,2)(0,4.5)
    \psline[linewidth=.7pt,linecolor=black,ArrowInside=->]{-}(0,4.5)(0,2)
    \psline[linewidth=.7pt,linecolor=black,ArrowInside=->]{-}(2.5,4.5)(0,2)
    \psline[linewidth=.9pt,linecolor=black]{->}(0,0)(0,-4.5)
    \psline[linewidth=.7pt,linecolor=black]{->}(0,0)(-4.5,-4.5)
    
    %\psplot[linecolor=black]{0}{3.3}{x}
    %\rput(*2.5 {x-0.5}){$\times\rho$}
    %\rput(*-0.05 {-999*x}){$\times\tau$}
    %\uput[0](-2.75,-1){$m_\rho$}
    %\uput[0](0,-1){$m_\tau$}
    %\uput[0](0,2){$m_\theta$}
%    \uput[0](1.6,3.5){$\times1/\rho$}
%    \uput[0](-2,3.5){$\times1/\tau$}
    \uput[0](-5,-3.5){$\times\rho$}
    \uput[0](-0.25,-3.5){$\times\tau$}
    \pscircle[fillstyle=solid,fillcolor=black,dimen=inner](0,0){0.07}

    \pscircle[fillstyle=solid,fillcolor=black,dimen=inner](0,2){0.07}
    \pscircle[fillstyle=solid,fillcolor=black,dimen=inner](0,3){0.07}
    \pscircle[fillstyle=solid,fillcolor=black,dimen=inner](1,3){0.07}

    \pscircle[fillstyle=solid,fillcolor=black,dimen=inner](0,4){0.07}
    \pscircle[fillstyle=solid,fillcolor=black,dimen=inner](1,3){0.07}
    \pscircle[fillstyle=solid,fillcolor=black,dimen=inner](1,4){0.07}
    \pscircle[fillstyle=solid,fillcolor=black,dimen=inner](2,4){0.07}
    
    \pscircle[fillstyle=solid,fillcolor=black,dimen=inner](0,-1){0.07}
    \pscircle[fillstyle=solid,fillcolor=black,dimen=inner](0,-2){0.07}
    \pscircle[fillstyle=solid,fillcolor=black,dimen=inner](0,-3){0.07}
    \pscircle[fillstyle=solid,fillcolor=black,dimen=inner](-1,-1){0.07}
    \pscircle[fillstyle=solid,fillcolor=black,dimen=inner](-1,-2){0.07}
    \pscircle[fillstyle=solid,fillcolor=black,dimen=inner](-1,-3){0.07}
    \pscircle[fillstyle=solid,fillcolor=black,dimen=inner](-2,-2){0.07}
    \pscircle[fillstyle=solid,fillcolor=black,dimen=inner](-2,-3){0.07}
    \pscircle[fillstyle=solid,fillcolor=black,dimen=inner](-3,-3){0.07}

\end{pspicture}
         \caption{\small $\Mop_2$.}
         \label{fig:M2*}
     \end{subfigure}
     \hfill
     \begin{subfigure}[t]{0.45\textwidth}
         \centering
\begin{pspicture}(-3,-3)(3,3)
\psset{unit = 0.5}
    \psaxes[linewidth=.4pt,linecolor=gray,showorigin=false,ticks=none,labels=none]{->}(0,0)(-5,-5)(5,5)[$p$,0][$q$,90]
%    \psline[linewidth=.7pt,linecolor=black]{->}(4,4)(0,0)
    \psline[linewidth=.7pt,linecolor=black,ArrowInside=->,arrowscale=1.5]{-}(4,4)(0,0)
    \psline[linewidth=.85pt,linecolor=black]{->}(0,4.5)(0,0)
    \psline[linewidth=.85pt,linecolor=black]{->}(0,0)(0,-4.5)
    \psline[linewidth=.7pt,linecolor=black]{->}(4,4.5)(4,0)
    \psline[linewidth=.7pt,linecolor=black]{->}(4,0)(4,-4.5)

    %\psplot[linecolor=black]{0}{3.3}{x}
    %\rput(*2.5 {x-0.5}){$\times\rho$}
    %\rput(*-0.05 {-999*x}){$\times\tau$}
%    \uput[0](0,0){$a_0$}
%    \uput[0](1,0){$a_1$}
    \uput[0](4,-0.5){$n$}
    
    \uput[0](1.75,4){$\times\rho$}
    \uput[0](-1.5,-3.5){$\times\tau$}
%    \uput[0](-2.25,3.5){$\times1/\tau$}
    
    \pscircle[fillstyle=solid,fillcolor=black,dimen=inner](0,0){0.07}
    \pscircle[fillstyle=solid,fillcolor=black,dimen=inner](1,0){0.07}
    \pscircle[fillstyle=solid,fillcolor=black,dimen=inner](2,0){0.07}
    \pscircle[fillstyle=solid,fillcolor=black,dimen=inner](3,0){0.07}
    \pscircle[fillstyle=solid,fillcolor=black,dimen=inner](4,0){0.07}
    
    \pscircle[fillstyle=solid,fillcolor=black,dimen=inner](0,1){0.07}
    \pscircle[fillstyle=solid,fillcolor=black,dimen=inner](1,1){0.07}
    \pscircle[fillstyle=solid,fillcolor=black,dimen=inner](2,1){0.07}
    \pscircle[fillstyle=solid,fillcolor=black,dimen=inner](3,1){0.07}
    \pscircle[fillstyle=solid,fillcolor=black,dimen=inner](4,1){0.07}
    \pscircle[fillstyle=solid,fillcolor=black,dimen=inner](0,2){0.07}
    \pscircle[fillstyle=solid,fillcolor=black,dimen=inner](1,2){0.07}
    \pscircle[fillstyle=solid,fillcolor=black,dimen=inner](2,2){0.07}
    \pscircle[fillstyle=solid,fillcolor=black,dimen=inner](3,2){0.07}
    \pscircle[fillstyle=solid,fillcolor=black,dimen=inner](4,2){0.07}
    \pscircle[fillstyle=solid,fillcolor=black,dimen=inner](0,3){0.07}
    \pscircle[fillstyle=solid,fillcolor=black,dimen=inner](1,3){0.07}
    \pscircle[fillstyle=solid,fillcolor=black,dimen=inner](2,3){0.07}
    \pscircle[fillstyle=solid,fillcolor=black,dimen=inner](3,3){0.07}
    \pscircle[fillstyle=solid,fillcolor=black,dimen=inner](4,3){0.07}
    \pscircle[fillstyle=solid,fillcolor=black,dimen=inner](0,-1){0.07}
    \pscircle[fillstyle=solid,fillcolor=black,dimen=inner](1,-1){0.07}
    \pscircle[fillstyle=solid,fillcolor=black,dimen=inner](2,-1){0.07}
    \pscircle[fillstyle=solid,fillcolor=black,dimen=inner](3,-1){0.07}
    \pscircle[fillstyle=solid,fillcolor=black,dimen=inner](4,-1){0.07}
    \pscircle[fillstyle=solid,fillcolor=black,dimen=inner](0,-2){0.07}
    \pscircle[fillstyle=solid,fillcolor=black,dimen=inner](1,-2){0.07}
    \pscircle[fillstyle=solid,fillcolor=black,dimen=inner](2,-2){0.07}
    \pscircle[fillstyle=solid,fillcolor=black,dimen=inner](3,-2){0.07}
    \pscircle[fillstyle=solid,fillcolor=black,dimen=inner](4,-2){0.07}
    \pscircle[fillstyle=solid,fillcolor=black,dimen=inner](0,-3){0.07}
    \pscircle[fillstyle=solid,fillcolor=black,dimen=inner](1,-3){0.07}
    \pscircle[fillstyle=solid,fillcolor=black,dimen=inner](2,-3){0.07}
    \pscircle[fillstyle=solid,fillcolor=black,dimen=inner](3,-3){0.07}
    \pscircle[fillstyle=solid,fillcolor=black,dimen=inner](4,-3){0.07}

\end{pspicture}
         \caption{\small $\Aop_n$.% (for further details, see Figure \ref{fig:An}).
         %The elements $a_0$ and $a_1$ refer to those mention in the proof of Lemma \ref{lem:duals}.
         }
         \label{fig:A4*}
     \end{subfigure}
     \hfill
\caption{\small Representations of $\Mop_2$ and $\Aop_n$ regarded as $\M_2$-modules.
     The vertical and oblique arrows indicate, respectively, multiplication by $\tau$ and by $\rho$.}
        \label{fig:duals}
\end{figure}
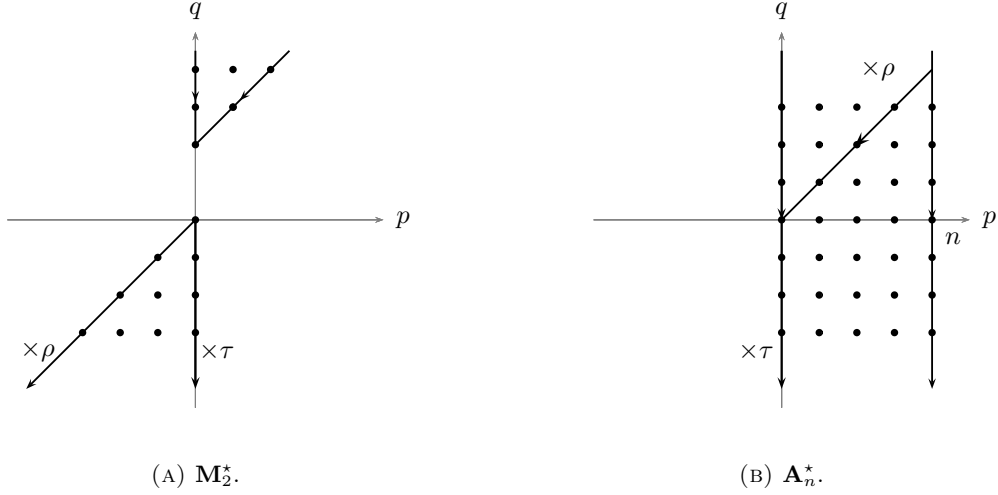
\end{ex}

\begin{rem}
\ \begin{enumerate}

\item The $\M_2$-module structure in bigraded Bredon homology is actually a graded module structure over $\M_2^\op$.     

\item The Poincaré duality isomorphism \eqref{eq:def_dual} is an isomorphism of bigraded $\M_2^\op$-modules as follows:
\begin{equation}
\label{eq:Poincare_duality_module_structure}
H^{C_2}_{*,*}(X;\uFt)\cong\Sigma^{2n,n}\left\{H^{*,*}_{C_2}(X;\uFt)^\op\right\}.
\end{equation}

\end{enumerate}
\end{rem}

From \eqref{eq:homology_structure}, \eqref{eq:Poincare_duality_module_structure} and Example~\ref{ex:M2_An_duals}, we can deduce a certain type of symmetry in the indices $p_i, q_i, r_j$ and $n_j$ in the decomposition \eqref{eq:structure_theorem}, which we state now along with a few other restrictions on this decomposition.

\begin{thm}
Let $X$ be a connected finite $C_2$-CW-complex that is a compact Real manifold of dimension $n$. Suppose that $X$ has a bigraded Bredon  cohomology decomposition as in \eqref{eq:structure_theorem}. Then the following conditions hold
\begin{enumerate}
\item For each $i\in I$, $q_i\leq n$;
\item For each $j\in J$, $r_j+n_j\leq 2n$;
\item If $X$ has a fixed point, then $r_j>0$  and $r_j+n_j<2n$ for all $j\in J$.
\end{enumerate}
Furthermore, the bigraded rank function $\rk^{p,q}(X)=\rk^{p,q}(H_{C_2}^{*,*}(X;\uFt))$ and the $a$-rank function $\rka^{p,q}(X):=\rka^{p,q}(H_{C_2}^{*,*}(X;\uFt))$  satisfy the following symmetry relations:

%Furthermore, the bigraded rank and $a$-rank functions  (see Definition~\ref{def:bigraded_ranks}) of %$N:=H_{C_2}^{*,*}(X;\uFt)$ satisfy the following symmetry relations:
%\begin{equation}
%k_{p,q} = |\{i\in I: (p_i,q_i)=(p,q) \}| \quad\text{and}\quad %m_{s,t}= |\{ j\in J: r_j=s \wedge n_j=t \} |,
%\label{eq:symmetry_restrictions}
%\end{equation}
%we have
\begin{equation}
\rk^{p,q}(X) = \rk^{2n-p,n-q}(X) 
\quad\text{ and }\quad 
\rka^{s,t}(X)=
\rka^{2n-s-t,t}(X), \qquad \forall p,q,s,t\in\Z_{\geq 0}.
\label{eq:restrictions_decomposition}
\end{equation}
\label{thm:restrictions_decomposition}
\end{thm}
\begin{proof}
The condition $r_j+n_j\leq 2n$ follows from \eqref{eq:singular_cohomology} by dimensional reasons. If $x_0\in X^{C_2}$ then the inclusion $i\colon x_0\to X$ gives a splitting
\[
H_{C_2}^{*,*}(X;\uFt) = \M_2\oplus \widetilde{H}_{C_2}^{*,*}(X;\uFt).
\]
Hence there exists an index $i\in I$ such that $(p_i,q_i)=(0,0).$ Since $H^0_{sing}(X;\F_2)=\F_2$, it now follows from 
\eqref{eq:singular_cohomology} that we must have $r_j>0,$ for all $j$.

By \eqref{eq:Poincare_duality_module_structure}, we have 
\begin{align*}
H^{C_2}_{*,*}(X;\uFt)&\cong{ \large \Sigma}^{2n,n}\left\{H^{*,*}_{C_2}(X;\uFt)^\op\right\}\\
&\cong  
\bigoplus_{i\in I} \Sigma^{2n-p_i,n-q_i}\M_2^\op \;\oplus\; \bigoplus_{j\in J} \Sigma^{2n-r_j,-n}\A_{n_j}^{\!\op} \\
&\cong  
\bigoplus_{i\in I} \Sigma^{2n-p_i,n-q_i}\M_2^\op \;\oplus\; (\bigoplus_{j\in J} \Sigma^{2n-r_j,0}\A_{n_j}^{\!\op}
\end{align*}
as $\M_2^\op$-modules. On the other hand, from \eqref{eq:homology_structure} and Example~\ref{ex:M2_An_duals}, we have
\begin{align*}
H^{C_2}_{*,*}(X;\uFt) &\cong \bigoplus_{i\in I} \Sigma^{p_i,q_i}\M_2^\star \;\oplus\; \bigoplus_{j\in J} \Sigma^{r_j,0}\A_{n_j}^{\!\star}\\
&\cong \bigoplus_{i\in I} \Sigma^{p_i,q_i}\M_2^\op \;\oplus\; \bigoplus_{j\in J} \Sigma^{r_j+n_j,0}\A_{n_j}^{\!\op}.
\end{align*}
Equating isomorphism types of $\M_2^\op$-modules in the two previous displays, we conclude that for each $i\in I$, there exists $i'\in I$ such that $(p_i,q_i)=(2n-p_{i'}, n-q_{i'})$; in particular $q_i\leq n$. Similarly, for each $j\in J$, there exists $j'\in J$ satisfying $n_{j'}=n_j$ and $r_{j'}+n_{j'}=2n-r_j$. By the definition of $\rk^{*,*}(N)$ and $\rka^{*,*}(N)$ in \eqref{eq:symmetry_restrictions}, the result now follows.
\end{proof}

%-----------------------------------------------------------------------------------------------------------------------------------------------
%
%	ILUSTRATION OF THE POINCARÉ SYMMETRIES IN BIGRADED BREDON COHOMOLOGY
%
%-----------------------------------------------------------------------------------------------------------------------------------------------

\subsection{Illustration of the Poincaré symmetries in bigraded Bredon cohomology}
\label{subsec:examples_bredon_cohomology}

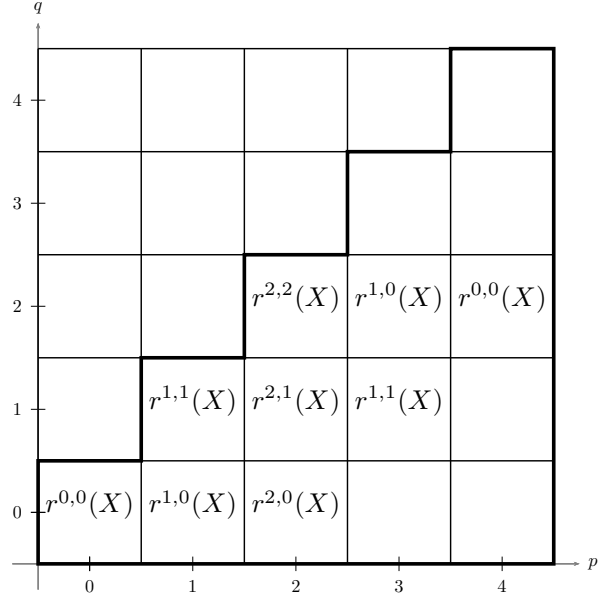
\begin{figure}[H]
\centering
\resizebox{0.45\textwidth}{!}{
\begin{pspicture}(-5,-6)(5,5.5)
\psset{unit = 1}

    \psaxes[linewidth=.4pt,linecolor=gray,showorigin=false,labels=none,ticks=none]{->}(-5,-5)(-5.5,-5.5)(5.5,5.5)[$p$,0][$q$,90]
    \psaxes[linewidth=0pt,linecolor=gray,dx=2,Ox=0]{-}(-4,-5)(-4,-5)(4,-5)
    \psaxes[linewidth=0pt,linecolor=gray,dy=2,Oy=0]{-}(-5,-4)(-5,-4)(-5,4)

%GRID
\psline{-}(-5,-5)(5,-5)
\psline{-}(-5,-3)(5,-3)
\psline{-}(-5,-1)(5,-1)
\psline{-}(-5,1)(5,1)
\psline{-}(-5,3)(5,3)
\psline{-}(-5,5)(5,5)
\psline{-}(-5,-5)(-5,5)
\psline{-}(-3,-5)(-3,5)
\psline{-}(-1,-5)(-1,5)
\psline{-}(1,-5)(1,5)
\psline{-}(3,-5)(3,5)
\psline{-}(5,-5)(5,5)

\pspolygon[linewidth=2pt](-5,-5)(5,-5)(5,5)(3,5)(3,3)(1,3)(1,1)(-1,1)(-1,-1)(-3,-1)(-3,-3)(-5,-3)

%LABELS
%zeroth row
\rput[B](-4,-4){\scalebox{1.5}{$r^{0,0}(X)$}}
\rput[B](-2,-4){\scalebox{1.5}{$r^{1,0}(X)$}}
\rput[B](0,-4){\scalebox{1.5}{$r^{2,0}(X)$}}
%\rput(! \pstnodescale\space 2 -4){$\bigtimes$}
%\rput(! \pstnodescale\space 4 -4){$\bigtimes$}
%first row
%\rput(! \pstnodescale\space -4 -2){$\bigtimes$}
\rput[B](-2,-2){\scalebox{1.5}{$\displaystyle r^{1,1}(X)$}}
\rput[B](0,-2){\scalebox{1.5}{$\displaystyle r^{2,1}(X)$}}
\rput[B](2,-2){\scalebox{1.5}{$\displaystyle r^{1,1}(X)$}}
%\rput(! \pstnodescale\space 4 -2){$\bigtimes$}
%second row
%\rput(! \pstnodescale\space -4 0){$\bigtimes$}
%\rput(! \pstnodescale\space -2 0){$\bigtimes$}
\rput[B](0,0){\scalebox{1.5}{$\displaystyle r^{2,2}(X)$}}
\rput[B](2,0){\scalebox{1.5}{$\displaystyle r^{1,0}(X)$}}
\rput[B](4,0){\scalebox{1.5}{$\displaystyle r^{0,0}(X)$}}
%third row
%\rput(! \pstnodescale\space 2 2){$\bigtimes$}
%\rput(! \pstnodescale\space 4 2){$\bigtimes$}
%fourth row
%\rput(! \pstnodescale\space 4 4){$\bigtimes$}

%\psline{-}(-4,-2)(-4,0)
%\psline{-}(-4,0)(-2,0)
%\psline{-}(2,-4)(4,-4)
%\psline{-}(4,-4)(4,-2)

%symmetry lines
%\psline[linestyle=dotted,dotsep=4pt]{-}(-5.5,-4.75)(5.5,0.75)
%\psline[linestyle=dotted,dotsep=3pt]{-}(-5.5,-5.3333)(5.5,1.333)
%\psline[linestyle=dotted,dotsep=4pt]{-}(-5.5,3.5)(3.5,-5.5)
    
\end{pspicture}
}
\caption{\small Symmetries of the bigraded rank function 
$r^{p,q}(X)=\rk^{p,q}H_{C_2}^{*,*}(X;\underline{\F_2})$
of a compact connected  Real manifold $X$  of dimension $n=2$ arising from Poincaré duality.
}
\label{fig:diamond}
\end{figure}

\begin{ex}[Poincaré duality in bigraded cohomology]
\label{ex:poincare_duality}
\ \begin{enumerate}%[ref=(\arabic*)]

\item
For a real smooth complete curve $X_g$ of genus $g$ such that $X_g(\C)$ has $r+1$ fixed curves (ovals), $r\leq g$ it is shown in 
\cite{Haz21} that
\[
H_{C_2}^{*,*}(X_g(\C);\uFt) = 
\M_2 \oplus (\Sigma^{1,0}\M_2)^r \oplus (\Sigma^{1,1}\M_2)^r \oplus (\Sigma^{1,0}\A_0)^{g-r}\oplus \Sigma^{2,1}\M_2.
\]    
\label{ex:riemann_surfaces}
Here the factor $(\Sigma^{1,1}\M_2)^r$ is dual to   $(\Sigma^{1,0}\M_2)^r$ in the sense that they correspond to each other under the symmetries of Theorem~\ref{thm:restrictions_decomposition} and $(\Sigma^{1,0}\A_0)^{g-r}$ is self-dual.

\item Denote by $\operatorname{SB}(2n+1)$ the real Severi--Brauer variety of dimension $2n+1$ without real points.
One can check that under the natural inclusion $\operatorname{SB}(1)\subset \operatorname{SB}(2n+1)$, there is a two-dimensional Real vector bundle $\operatorname{SB}(2n+1)\setminus\operatorname{SB}(1)\to \operatorname{SB}(2n-1)$.
This can be used to compute the bigraded cohomology inductively:
\[
H_{C_2}^{*,*}(\operatorname{SB}(2n+1);\uFt)
= \A_2\oplus \Sigma^{4,0}\A_2\oplus\dotsm\oplus\Sigma^{4n-4,0}\A_2
\oplus\Sigma^{4n,0}\A_2
\]
(see \cite{pedro&paulo:quaternionic_algebraic_cycles}) and we see that  each $\Sigma^{4k}\A_2$ factor has $\Sigma^{4(n-k)}\A_2$ as dual.

\item The twisted projective plane, sometimes denoted by $\RP_{tw}^2$, is the real projective plane $\RP^2$ endowed with the $C_2$-action given by a rotation of 180 degrees. %(see Figure~\ref{fig:twisted_plane}).
It is not a Real manifold, as  its bigraded cohomology (see \cite[Example 3.2]{clover_may:freeness_theorem_cohomology})
\[
H_{C_2}^{*,*}(\RP_{tw}^2;\underline{\mathbf{F}_2}) = \M_2 \oplus \Sigma^{1,1}\M_2 \oplus\Sigma^{2,1}\M_2
\]
has no \(\Sigma^{1,0}\M_2\) factor corresponding to \(\Sigma^{1,1}\M_2.\)
\label{counterex:twisted_plane}

\begin{comment}
\begin{figure}[H]
     \centering
\resizebox{0.3\textwidth}{!}{
\begin{pspicture}(-3.5,-3)(3.5,3)
\psset{unit = 0.5}
    \psaxes[linewidth=.4pt,linecolor=gray,showorigin=false,ticks=none,labels=none]{->}(0,0)(-5.5,-5.5)(5.5,5.5)[$x$,0][$y$,90]
    
    \pscircle(0,0){5}
    
\psarc[linewidth=0pt,arrowsize=6pt]{->}{5}{269}{271}
\psarc[linewidth=0pt,arrowsize=6pt]{->}{5}{89}{91}

\psarc[arrowsize=3pt]{->}{1.25}{-150}{150}
\rput[bl]{0}(1.25,1.25){$C_2$}

%\pscircle[linestyle=none,fillstyle=solid,fillcolor=blue](0,0){0.1}
\pscircle[fillstyle=solid,fillcolor=blue,dimen=inner](0,0){0.15}
\pscircle[fillstyle=solid,fillcolor=blue,dimen=inner](-5,0){0.15} 
\pscircle[fillstyle=solid,fillcolor=blue,dimen=inner](5,0){0.15}

\end{pspicture}
}
         \caption{\small Twisted projective plane $\RP_{tw}^2$ of Example \ref{ex:poincare_duality} (\ref{counterex:twisted_plane}) (see \cite{Haz21} and \cite{clover_may:freeness_theorem_cohomology} for more details).}
         \label{fig:twisted_plane}
\end{figure}
\end{comment}

\end{enumerate}
\end{ex}

%\textcolor{red}{We need to reference Examples \ref{ex:riemann_surfaces} and \ref{ex:brauer-severi}.}

%-----------------------------------------------------------------------------------------------------------------------------------------------
%
%	APPLICATIONS OF THE POINCARÉ SYMMETRIES TO THE GM-CONDITION
%
%-----------------------------------------------------------------------------------------------------------------------------------------------
\subsection{Applications of the Poincaré symmetries to the GM-condition}
\label{subsec:applications_poincare_symmetries}

We can rederive the following result of Krasnov as a consequence of the restrictions on decomposition \eqref{eq:structure_theorem} imposed by Theorem \ref{thm:restrictions_decomposition}.

\begin{prop}[Krasnov \cite{krasnov:harnack-thom_inequalities}]
\label{prop:krasnov}
Let $V$ be a non-singular complete real algebraic surface such that $V(\R)$ is non-empty and $H_{sing}^1(V(\C);\F_2)$ is trivial.
Then $V$ is a GM-variety.
\end{prop}
\begin{proof}
Since $H_{sing}^1(V(\C);\F_2)=0$ and $V(\R)\neq\varnothing$,  
condition {\it (3)} in Theorem~\ref{thm:restrictions_decomposition} applies and hence the list of allowed pairs $(r_j, n_j)$ in decomposition in \eqref{eq:structure_theorem}
is:
\[
L=\{(1,0), (2,0), (3,0), (1,1), (2,1), (1,2)\}.
\]

If $(2,1)\in L$ or $(1,2)\in L$, we would have $H_{sing}^3(V(\C);\F_2)\neq 0$, contradicting  the triviality of $H_{sing}^1(V(\C);\F_2)$, by Poincar\'{e} duality.

Finally, if $(1,1)\in L$, then, by the second symmetry condition in \eqref{eq:symmetry_restrictions}, we would have $(2,1)\in L$, which is not possible.
\end{proof}

\begin{ex}[K3 surface whose real part is \(S^2\amalg S^2\)]
\label{ex:K3_surfaces}
A K3 surface is a complete complex algebraic surface with trivial canonical class and with first Betti number equal to zero. All K3 surfaces have the same Hodge polynomial, given by
\begin{equation}
\label{eq:Hodge-K3}
H(u,v) = 1 + u^2 + v^2 + 20uv + u^2v^2.
\end{equation}
From Proposition \ref{prop:krasnov}, it follows that all real K3 surfaces are Galois-maximal varieties.

As shown in Silhol \cite{silhol:real_algebraic_surfaces}, real K3 surfaces are completely determined, up to isomorphism, by the topology of the real part.
More precisely, their moduli space is described by values of the total Betti number $b_*(S(\R))$ and the Euler characteristic $\chi(S(\R))$ of the real points, both computed with $\F_2$ coefficients (see p.~189 on \textit{loc.~cit.}).

We use the results of Theorem \ref{thm:cohomology_fixed_points} and Corollary \ref{cor:singular_cohomology} together with Theorems \ref{thm:main_result} and \ref{thm:restrictions_decomposition} to compute the bigraded cohomology of the K3 surface whose real part consists of the union of two spheres.
As we can see in Figure \ref{fig:K3_surfaces} below, this case corresponds to the values $b_*(S(\R))=\chi(S(\R))=4$, and we have:
\begin{equation}
\begin{cases}
\dim H_{sing}^0(S(\R))=\dim H_{sing}^2(S(\R)) = 2 \\
\dim H_{sing}^1(S(\R))=0.
\end{cases}
\label{eq:dimensionK3}
\end{equation}

We will start by computing the free part.
As there is one real point, Theorem \ref{thm:main_result} guarantees the existence of the factors $\M_2$ and its Poincaré dual $\Sigma^{4,2}\M_2$ in the bigraded cohomology of $S(\C)$.
This is not enough to generate completely the singular cohomology of $S(\R)$.
If we look at \eqref{eq:dimensionK3}, the results of Theorem \ref{thm:cohomology_fixed_points} indicate the existence of other two factors $\Sigma^{p_1,q_1}\M_2$ and $\Sigma^{p_2,q_2}\M_2$ such that $p_1=q_1$ and $p_2-q_2=2$.
In addition, equivariant Poincaré duality implies that $p_1+p_2=4$ and $q_1+q_2=2$.
Finally, the fact $b_1(S(\C))=0$ implies that $p_i\neq 1,3$.
Hence, the only possibilities are 
\[
(p_1,q_1)=(2,0) \quad \mbox{ and } \quad (p_2,q_2)=(2,2),
\]
and the free part is completely described.

The non-free part, is given by Corollary \ref{cor:singular_cohomology} as follows: From \eqref{eq:Hodge-K3} we can see that all $b_2=22$ for all K3 surfaces; as each free orbit contributes with a factor two, there must be 10 different copies, resulting
\[
H_{C_2}^{*,*}(S(\C);\underline{\F_2}) = \M_2 \oplus \Sigma^{2,0}\M_2 \oplus\Sigma^{2,2}\M_2  \oplus \left(\Sigma^{2,0}\A_0\right)^{10} \oplus \Sigma^{4,2}\M_2.
\]
\end{ex}

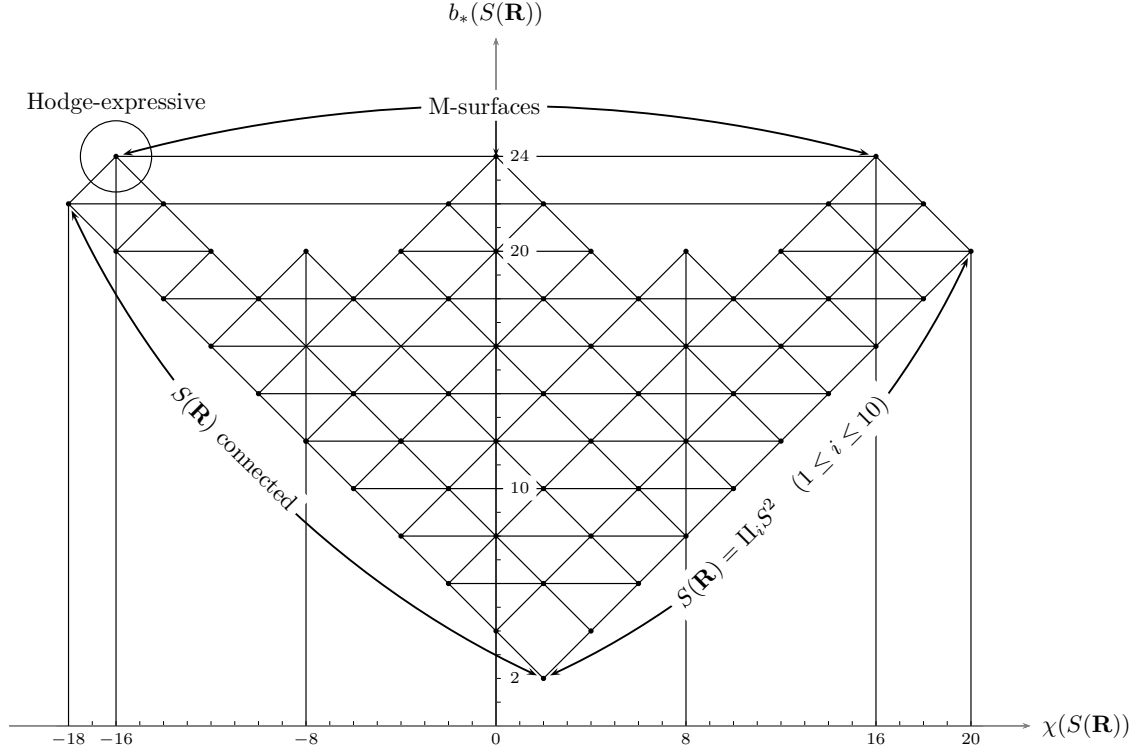
\begin{figure}[H]
     \centering
\resizebox{0.65\textwidth}{!}{
\begin{pspicture}(-5.5,-0.5)(5.5,10.75)
\psset{unit = 0.35}
\psaxes[linewidth=0pt,linecolor=gray,ticksize=2pt,showorigin=false,labels=none]{-}(0,0)(-18.5,0)(20.5,24.5)
\psaxes[linewidth=0.4pt,linecolor=gray,ticks=none,showorigin=false,labels=none]{->}(0,0)(-20.5,0)(22.5,29)[$\chi(S(\R))$,0][$b_*(S(\R))$,90]

%2's tier
\pscircle[fillstyle=solid,fillcolor=black,dimen=inner](2,2){0.07}
%4's tier
\pscircle[fillstyle=solid,fillcolor=black,dimen=inner](4,4){0.07}
\pscircle[fillstyle=solid,fillcolor=black,dimen=inner](0,4){0.07}
%6's tier
\pscircle[fillstyle=solid,fillcolor=black,dimen=inner](6,6){0.07}
\pscircle[fillstyle=solid,fillcolor=black,dimen=inner](2,6){0.07}
\pscircle[fillstyle=solid,fillcolor=black,dimen=inner](-2,6){0.07}
%8's tier
\pscircle[fillstyle=solid,fillcolor=black,dimen=inner](8,8){0.07}
\pscircle[fillstyle=solid,fillcolor=black,dimen=inner](4,8){0.07}
\pscircle[fillstyle=solid,fillcolor=black,dimen=inner](0,8){0.07}
\pscircle[fillstyle=solid,fillcolor=black,dimen=inner](-4,8){0.07}
%10's tier
\pscircle[fillstyle=solid,fillcolor=black,dimen=inner](10,10){0.07}
\pscircle[fillstyle=solid,fillcolor=black,dimen=inner](6,10){0.07}
\pscircle[fillstyle=solid,fillcolor=black,dimen=inner](2,10){0.07}
\pscircle[fillstyle=solid,fillcolor=black,dimen=inner](-2,10){0.07}
\pscircle[fillstyle=solid,fillcolor=black,dimen=inner](-6,10){0.07}
%12's tier
\pscircle[fillstyle=solid,fillcolor=black,dimen=inner](12,12){0.07}
\pscircle[fillstyle=solid,fillcolor=black,dimen=inner](8,12){0.07}
\pscircle[fillstyle=solid,fillcolor=black,dimen=inner](4,12){0.07}
\pscircle[fillstyle=solid,fillcolor=black,dimen=inner](0,12){0.07}
\pscircle[fillstyle=solid,fillcolor=black,dimen=inner](-4,12){0.07}
\pscircle[fillstyle=solid,fillcolor=black,dimen=inner](-8,12){0.07}
%14's tier
\pscircle[fillstyle=solid,fillcolor=black,dimen=inner](14,14){0.07}
\pscircle[fillstyle=solid,fillcolor=black,dimen=inner](10,14){0.07}
\pscircle[fillstyle=solid,fillcolor=black,dimen=inner](6,14){0.07}
\pscircle[fillstyle=solid,fillcolor=black,dimen=inner](2,14){0.07}
\pscircle[fillstyle=solid,fillcolor=black,dimen=inner](-2,14){0.07}
\pscircle[fillstyle=solid,fillcolor=black,dimen=inner](-6,14){0.07}
\pscircle[fillstyle=solid,fillcolor=black,dimen=inner](-10,14){0.07}
%16's tier
\pscircle[fillstyle=solid,fillcolor=black,dimen=inner](16,16){0.07}
\pscircle[fillstyle=solid,fillcolor=black,dimen=inner](12,16){0.07}
\pscircle[fillstyle=solid,fillcolor=black,dimen=inner](8,16){0.07}
\pscircle[fillstyle=solid,fillcolor=black,dimen=inner](4,16){0.07}
\pscircle[fillstyle=solid,fillcolor=black,dimen=inner](0,16){0.07}
\pscircle[fillstyle=solid,fillcolor=black,dimen=inner](-4,12){0.07}
\pscircle[fillstyle=solid,fillcolor=black,dimen=inner](-8,12){0.07}
\pscircle[fillstyle=solid,fillcolor=black,dimen=inner](-12,16){0.07}
%18's tier
\pscircle[fillstyle=solid,fillcolor=black,dimen=inner](18,18){0.07}
\pscircle[fillstyle=solid,fillcolor=black,dimen=inner](14,18){0.07}
\pscircle[fillstyle=solid,fillcolor=black,dimen=inner](10,18){0.07}
\pscircle[fillstyle=solid,fillcolor=black,dimen=inner](6,18){0.07}
\pscircle[fillstyle=solid,fillcolor=black,dimen=inner](2,18){0.07}
\pscircle[fillstyle=solid,fillcolor=black,dimen=inner](-2,18){0.07}
\pscircle[fillstyle=solid,fillcolor=black,dimen=inner](-6,18){0.07}
\pscircle[fillstyle=solid,fillcolor=black,dimen=inner](-10,18){0.07}
\pscircle[fillstyle=solid,fillcolor=black,dimen=inner](-14,18){0.07}
%20's tier
\pscircle[fillstyle=solid,fillcolor=black,dimen=inner](20,20){0.07}
\pscircle[fillstyle=solid,fillcolor=black,dimen=inner](16,20){0.07}
\pscircle[fillstyle=solid,fillcolor=black,dimen=inner](12,20){0.07}
\pscircle[fillstyle=solid,fillcolor=black,dimen=inner](8,20){0.07}
\pscircle[fillstyle=solid,fillcolor=black,dimen=inner](4,20){0.07}
\pscircle[fillstyle=solid,fillcolor=black,dimen=inner](0,20){0.07}
\pscircle[fillstyle=solid,fillcolor=black,dimen=inner](-4,20){0.07}
\pscircle[fillstyle=solid,fillcolor=black,dimen=inner](-8,20){0.07}
\pscircle[fillstyle=solid,fillcolor=black,dimen=inner](-12,20){0.07}
\pscircle[fillstyle=solid,fillcolor=black,dimen=inner](-16,20){0.07}
%22's tier
\pscircle[fillstyle=solid,fillcolor=black,dimen=inner](2,22){0.07}
\pscircle[fillstyle=solid,fillcolor=black,dimen=inner](-2,22){0.07}
\pscircle[fillstyle=solid,fillcolor=black,dimen=inner](-14,22){0.07}
\pscircle[fillstyle=solid,fillcolor=black,dimen=inner](-18,22){0.07}
\pscircle[fillstyle=solid,fillcolor=black,dimen=inner](18,22){0.07}
\pscircle[fillstyle=solid,fillcolor=black,dimen=inner](14,22){0.07}
%24's tier
\pscircle[fillstyle=solid,fillcolor=black,dimen=inner](0,24){0.07}
\pscircle[fillstyle=solid,fillcolor=black,dimen=inner](16,24){0.07}
\pscircle[fillstyle=solid,fillcolor=black,dimen=inner](-16,24){0.07}

%diagonal lines
\psline[linewidth=.5pt,linecolor=black]{-}(2,2)(20,20)
\psline[linewidth=.5pt,linecolor=black]{-}(0,4)(18,22)
\psline[linewidth=.5pt,linecolor=black]{-}(-2,6)(16,24)
\psline[linewidth=.5pt,linecolor=black]{-}(-4,8)(8,20)
\psline[linewidth=.5pt,linecolor=black]{-}(-6,10)(4,20)
\psline[linewidth=.5pt,linecolor=black]{-}(-8,12)(2,22)
\psline[linewidth=.5pt,linecolor=black]{-}(-10,14)(0,24)
\psline[linewidth=.5pt,linecolor=black]{-}(-12,16)(-8,20)
\psline[linewidth=.5pt,linecolor=black]{-}(-14,18)(-12,20)
\psline[linewidth=.5pt,linecolor=black]{-}(-16,20)(-14,22)
\psline[linewidth=.5pt,linecolor=black]{-}(-18,22)(-16,24)
%horizontal lines
\psline[linewidth=.5pt,linecolor=black]{-}(-2,6)(6,6)
\psline[linewidth=.5pt,linecolor=black]{-}(-4,8)(8,8)
\psline[linewidth=.5pt,linecolor=black]{-}(-6,10)(10,10)
\psline[linewidth=.5pt,linecolor=black]{-}(-8,12)(12,12)
\psline[linewidth=.5pt,linecolor=black]{-}(-10,14)(14,14)
\psline[linewidth=.5pt,linecolor=black]{-}(-12,16)(16,16)
\psline[linewidth=.5pt,linecolor=black]{-}(-14,18)(18,18)
\psline[linewidth=.5pt,linecolor=black]{-}(-4,20)(4,20)
\psline[linewidth=.5pt,linecolor=black]{-}(-16,20)(-12,20)
\psline[linewidth=.5pt,linecolor=black]{-}(12,20)(20,20)
\psline[linewidth=.5pt,linecolor=black]{-}(-18,22)(18,22)
\psline[linewidth=.5pt,linecolor=black]{-}(-16,24)(16,24)
%anti-diagonal lines
\psline[linewidth=.5pt,linecolor=black]{-}(-18,22)(2,2)
\psline[linewidth=.5pt,linecolor=black]{-}(-16,24)(4,4)
\psline[linewidth=.5pt,linecolor=black]{-}(-8,20)(6,6)
\psline[linewidth=.5pt,linecolor=black]{-}(-4,20)(8,8)
\psline[linewidth=.5pt,linecolor=black]{-}(-2,22)(10,10)
\psline[linewidth=.5pt,linecolor=black]{-}(0,24)(12,12)
\psline[linewidth=.5pt,linecolor=black]{-}(8,20)(14,14)
\psline[linewidth=.5pt,linecolor=black]{-}(12,20)(16,16)
\psline[linewidth=.5pt,linecolor=black]{-}(14,22)(18,18)
\psline[linewidth=.5pt,linecolor=black]{-}(16,24)(20,20)
%vertical lines
\psline[linewidth=.5pt,linecolor=black]{-}(-18,0)(-18,22)
\psline[linewidth=.5pt,linecolor=black]{-}(-16,0)(-16,24)
\psline[linewidth=.5pt,linecolor=black]{-}(-8,0)(-8,20)
\psline[linewidth=.5pt,linecolor=black]{-}(0,0)(0,24)
\psline[linewidth=.5pt,linecolor=black]{-}(8,0)(8,20)
\psline[linewidth=.5pt,linecolor=black]{-}(16,0)(16,24)
\psline[linewidth=.5pt,linecolor=black]{-}(20,0)(20,20)

%ticks on x-axis
\rput[t]{0}(0,-0.3){$\scriptstyle 0$}
\rput[t]{0}(-8,-0.3){$\scriptstyle -8$}
\rput[t]{0}(-16,-0.3){$\scriptstyle -16$}
\rput[t]{0}(-18,-0.3){$\scriptstyle -18$}
\rput[t]{0}(8,-0.3){$\scriptstyle 8$}
\rput[t]{0}(16,-0.3){$\scriptstyle 16$}
\rput[t]{0}(20,-0.3){$\scriptstyle 20$}
%ticks on y-axis
\rput[Bl]{0}(0.3,1.8){\psframebox*{$\scriptstyle 2$}}
\rput[Bl]{0}(0.3,9.8){\psframebox*{$\scriptstyle 10$}}
\rput[Bl]{0}(0.3,19.8){\psframebox*{$\scriptstyle 20$}}
\rput[Bl]{0}(0.3,23.8){\psframebox*{$\scriptstyle 24$}}

\psline[linewidth=.5pt,linecolor=black]{->}(0,26)(0,24)
\pnode(2,2){A}
\pnode(20,20){B}
\pnode(-18,22){C}
\pnode(16,24){D}
\pnode(-16,24){E}
\psset{nodesep=3pt}
\ncarc[arcangle=20]{<->}{B}{A}
\ncarc[arcangle=20]{<->}{A}{C}
\ncarc[arcangle=15]{<->}{E}{D}
%\psframebox*{\psline[linewidth=0.5pt,linearc=25]{<->}(2,2)(11,11)(20,20)}

\rput[B]{45}(12.25,9.25){\psframebox*{$S(\R)=\amalg_i S^2\quad (1\leq i\leq 10)$}}
\rput[t]{315}(-10.5,12.25){\psframebox*{$S(\R)$ connected}}
\rput[b]{0}(-0.5,25.5){\psframebox*{M-surfaces}}
\rput[b]{0}(-16,25.5){\psframebox*{Hodge-expressive}}
\pscircle[linewidth=0.5pt](-16,24){1.5}

\end{pspicture}
}

         \caption{\small Moduli space of real K3 surfaces as presented in \cite[VIII(3.3)]{silhol:real_algebraic_surfaces}.
         We have added the location of the Hodge-expressive surfaces.}
         \label{fig:K3_surfaces}
\end{figure}

\begin{ex}[General real K3 surfaces]
%\label{ex:}
A similar reasoning allows to compute the bigraded cohomology of a general real K3 surface $S$.
Setting
\begin{comment}
In general
\begin{equation*}
\begin{cases}
\displaystyle \dim H^0(S(\R))=\dim H^2(S(\R)) = \displaystyle \frac{b_*(S(\R))+\chi(S(\R))}{4} \\
\displaystyle \dim H^1(S(\R))= \displaystyle \frac{b_*(S(\R))-\chi(S(\R))}{2}.
\end{cases}
\label{eq:betti_numbers_k3_surfaces}
\end{equation*}
%
Therefore, by setting
%
\[
a=\frac{b_*(S(\R))+\chi(S(\R))}{4}-1=b_0(S(\R))-1, \quad b = \frac{b_*(S(\R))-\chi(S(\R))}{2}=b_1(S(\R)) \quad \mbox{and} \quad c= \frac{24-b_*(S(\R))}{2},
\]
%
and using the results of Theorem \ref{thm:}
\end{comment}
%
\[
a = b_0(S(\R))- 1, \quad b = b_1(S(\R)) \quad \mbox{and} \quad c=\frac{b_*(S(\C)) - b_*(S(\R))}{2},
\]
we obtain
\begin{equation}
H_{C_2}^{*,*}(S(\C);\underline{\F_2}) = \M_2  \oplus \left(\Sigma^{2,0}\M_2\right)^{a} \oplus \left(\Sigma^{2,1}\M_2 \right)^{b} \oplus \left(\Sigma^{2,2}\M_2\right)^{a} \oplus \left(\Sigma^{2,0}\A_0\right)^{c} \oplus \Sigma^{4,2}\M_2.
\label{eq:cohomology_K3-surfaces}
\end{equation}
\end{ex}

We complete our analysis of the bigraded Bredon cohomology of K3 surfaces by considering   the case of Hodge-expressive K3 surfaces.

\begin{ex}[Hodge-expressive K3 surfaces]
Recall that this family of varieties is defined by the absence of torsion in integral cohomology and the relation
$P_{S(\R)}(t)=H_{S(\C)}(t,1)$.
%The work of \cite{BruSch21} includes more information about those K3 surfaces which are Hodge-expressive.

The condition $P_{S(\R)}(t)=H_{S(\C)}(t,1)$ yields
\[
P_{S(\R)}(t)=2+20t+2t^2,
\]
and so Hodge-expressive K3 surfaces are completely characterized by the conditions
\[
b_*(S(\R))=24 \quad \mbox{and } \quad \chi(S(\R))=-16 .
\]
In particular, the equality $b_*(S(\R))=24$ says that the Hodge-expressive K3 surfaces are maximal, as required.

From \eqref{eq:cohomology_K3-surfaces}, the bigraded cohomology of this particular class of K3 surfaces is 
\begin{equation}
H_{C_2}^{*,*}(S(\C);\underline{\F_2}) = \M_2 \oplus \Sigma^{2,0}\M_2 \oplus \left(\Sigma^{2,1}\M_2\right)^{20} \oplus \Sigma^{2,2}\M_2 \oplus \Sigma^{4,2}\M_2 .
\label{eq:bredon_cohomology_HE_K3}
\end{equation}
Recall now the bigraded rank functions:
\[
r^{p,q} (S(\C)) = \rk^{p,q}(H_{C_2}^{*,*}(S(\C);\underline{\F_2}).
\]
Using \eqref{eq:bredon_cohomology_HE_K3}, we see that the Hodge-expressive K3 surfaces satisfy the especial relation
\begin{equation}
r^{p+q,q} (S(\C)) = h^{p,q} (S(\C))
\label{eq:hodge-expressive_bredon}
\end{equation}
for all $p,q\geq 0$. As mentioned in Remark~\ref{rem:Hodge-birank}, this relation always implies the Hodge-expressive condition.
%It turns out that this is a sufficient condition in general, \emph{i.e.}~if a Real M-variety $V$ has no torsion and satisfies \eqref{eq:hodge-expressive_bredon}, then $V$ is a Hodge-expressive variety.
%Indeed, the results of Theorem \ref{thm:cohomology_fixed_points} show that
%
%\[
%b_n(V(\R)) = \sum_{\substack{i\in I : \\ p_i-q_i=n}} k^{p_i,q_i};
%\]
%
%then the assertion follows from \eqref{eq:hodge_expressive_bis}.
\end{ex}

We finish this section by formulating a result
analogous to Proposition~\ref{prop:krasnov}
for real threefolds.

\begin{prop}
\label{prop:3fold_GM}
Let $V$ be a smooth irreducible real 3-fold such that:
\begin{enumerate}[(i)]
    \item $V(\R)\neq\varnothing$;
    \item the forgetful map $\psi^{4,2}\colon H_{C_2}^{4,2}(V(\C);\uFt)\to H_{sing}^{4}(V(\C);\F_2)$ is onto;
    \item $H^1_{sing}(V(\C);\F_2)=0$.
\end{enumerate}
Then $V$ is a GM-variety.
\end{prop}
\begin{proof} Let $L$ denote the list pairs $(r_j,n_j)$ appearing in the decomposition of \eqref{eq:structure_theorem} for the bigraded cohomology of $V(\C).$ That is, $L\subset\Z_{\geq 0}^2$ is the support of 
$\rka\bigl(H_{C_2}^{*,*}(V(\C);\uFt)\bigr).$
By Theorem~\ref{thm:main_result}, we only need to show that $L\subset\Z_{\geq 0}\times\{0\}.$ 

Since $H_{sing}^i(V(\C);\F_2)=0$ for $i=1,5$, it follows from \eqref{eq:singular_cohomology} that 
\[
(r,s)\in L \Rightarrow r, r+s \neq 1,5.
\]
Also, as $V(\R)\neq\varnothing$, condition (3) in Theorem~\ref{thm:restrictions_decomposition} applies. We conclude that
\[
L\subset\{(2,0),(3,0),(4,0),(2,1),(3,1),
(2,2)\}.
\]
We also have:
\begin{itemize}
    \item If  $(2,1)\in L$ then $(3,1)\in L$ by \eqref{eq:restrictions_decomposition} and $\psi^{4,2}$ is not onto;
    \item If $(3,1)\in L$ then $\psi^{4,2}$ is not onto;
    \item If $(2,2)\in L$ then  $\psi^{4,2}$ is not onto.
\end{itemize}
The result follows.
\end{proof}

\begin{ex}[Real threefolds]
\label{ex:threefolds}
Let $Y\subset\mathbf{P}^4$ be a smooth real threefold hypersurface containing a real line $L$, such as a real cubic (\emph{cf.}~\cite{DeMan00}).
Denote by $i$ and $j$ the inclusions $i\colon Y(\C)\hookrightarrow \mathbf{P}^4(\C)$ and
$j\colon L(\C)\hookrightarrow Y(\C),$ respectively, and by $h$ the generator of the summand $\Sigma^{2,1}\M_2$ in the decomposition
\eqref{eq:decomp_Pn} for $n=4$.
Since $(i\circ j)^*h$ is the generator of the summand $\Sigma^{2,1}\M_2$ in the decomposition
\eqref{eq:decomp_Pn} for $n=1$, it follows that $i^*h$ generates a 
summand of the form $\Sigma^{2,1}\M_2$ in the decomposition of the bigraded cohomology of $Y(\C).$ By Theorem~\ref{thm:restrictions_decomposition}, this decomposition has also a summand of the form $\Sigma^{4,2}\M_2.$ Now, by the Lefschtez hyperplane theorem,  $H^1_{sing}(Y(\C);\F_2)=0$ and
$H^4_{sing}(Y(\C);\F_2)=\F_2$. It follows that 
\[
\psi^{4,2}\colon H_{C_2}^{4,2}(Y(\C);\uFt)\to H_{sing}^{4}(Y(\C);\F_2)
\]
is onto and Proposition~\ref{prop:3fold_GM} applies. Hence, $Y$ is a GM-variety.
\end{ex}

\begin{ex}[Real cubic]

As a particular example of the above we consider the cubic $Y$ whose real part consists of the disjoint union of $S^3$ and $\RP^3$ (\emph{cf.}~\cite{krasnov:classification_real_cubics}).
All the real cubics have the same Hodge polynomial:
\begin{equation}
H(u,v) = 1 + uv + 5u^2v + 5uv^2 + u^2v^2 + u^3+v^3 .
\label{eq:hodge_polynomial_cubic}
\end{equation}

To compute the bigraded cohomology, we proceed as in Example \ref{ex:K3_surfaces} and find out first about the free part.
There must be factors $\M_2$ and $\Sigma^{6,3}\M_2$, and using the results of Example \ref{ex:threefolds}, we also deduce the existence of $\Sigma^{2,1}\M_2$ and $\Sigma^{4,2}\M_2$.
However, we still need two more factors to recover the cohomology of the real part, as
\begin{equation*}
\begin{cases}
\dim H^0(Y(\R) = \dim H^3(Y(\R)) = 2, \\
\dim H^1(Y(\R)) = \dim H^2(Y(\R)) = 1.
\end{cases}
\end{equation*}
Due to the restrictions imposed by the singular cohomology of $Y(\C)$ reflected in the Hodge polynomial of \eqref{eq:hodge_polynomial_cubic} and by equivariant Poincaré duality, the two remaining factors must be concentrated in topological dimension three, which provides $\Sigma^{3,0}\M_2$ and its Poincaré dual $\Sigma^{3,3}\M_2$.

Finally, to find the non-free part, we use the fact that $\dim H^3(Y(\C);\F_2)=10$:
\[
H_{C_2}^{*,*}(Y(\C);\uFt) = \M_2 \oplus \Sigma^{2,1}\M_2\oplus \Sigma^{3,0}\M_2 \oplus  \Sigma^{3,3}\M_2 \oplus \Sigma^{4,2}\M_2\oplus
\left(\Sigma^{3,0}\A_0\right)^4 \oplus \Sigma^{6,3}\M_2.
\]
\end{ex}

\bibliographystyle{acm}
\bibliography{data_base}
%----------------------------------------------------------------------------------------
%%%	   NO CITED	%%%
%\nocite{choquetmanifolds}
%----------------------------------------------------------------------------------------

\end{document}